\documentclass[11pt]{amsart}

\usepackage[foot]{amsaddr}
\usepackage{hyperref}
\usepackage{macros-arxiv}
\usepackage{cleveref}
\hypersetup{colorlinks=true,allcolors=darkblue}
\usepackage{comment} 
\usepackage{enumerate}
\usepackage[shortlabels]{enumitem}
\usepackage{overpic}
\usebiggeometry

\title[Accelerated, Optimal, and Parallel]{Accelerated, Optimal, and Parallel: \\
  Some results on model-based stochastic optimization}
\author{Karan Chadha$^{1, *}$}
\address[$*$]{Denotes equal contribution; authors listed in alphabetical order.}
\address[$1$]{Electrical Engineering Department, Stanford University, Stanford, CA}
\email{knchadha@stanford.edu}
\author{Gary Cheng$^{1, *}$}
\email{chenggar@stanford.edu}
\author{John C. Duchi$^{1,2}$}
\address[$2$]{Statistics Department, Stanford University, Stanford, CA}
\email{jduchi@stanford.edu}

\begin{document}
\maketitle

\begin{abstract}
  We extend the Approximate-Proximal Point (\aprox) family of model-based
  methods for solving stochastic convex optimization problems, including
  stochastic subgradient, proximal point, and bundle methods, to the
  minibatch and accelerated setting. To do so, we propose specific
  model-based
  algorithms and an acceleration scheme for which we provide non-asymptotic
  convergence guarantees, which are order-optimal in all
  problem-dependent constants and provide linear speedup in
  minibatch size, while maintaining the desirable robustness traits
  (e.g.\ to stepsize) of the \aprox family.
  Additionally, we show improved convergence rates and
  matching  lower bounds identifying new fundamental constants
  for ``interpolation'' problems, whose 
  importance in statistical machine learning is growing; this, for
  example, gives a parallelization strategy for alternating
  projections. We corroborate our theoretical results with
  empirical testing to demonstrate the gains accurate
  modeling, acceleration, and minibatching provide.
\end{abstract}

\section{Introduction}
\label{sec:intro}

We move beyond stochastic and ``minibatch''-gradient methods for
stochastic optimization problems to develop parallelizable
and minibatch aware model-based and
(approximate) proximal point methods for the problem
\begin{equation}
  \label{eqn:objective}
  \begin{split}
    \minimize ~ & f(x)
    \defeq \E_P[F(x;\statrv)]
    = \int_\statdomain F(x; \statval) dP(\statval) \\
    \subjectto ~ & x \in \domain.
  \end{split}
\end{equation}
Here, $\statdomain$ denotes the sample space, and $\statrv \sim P$ is an
$\statdomain$-valued random variable, where for each sample $\statval \in
\statdomain$, $F(\cdot; \statval): \R^n \rightarrow \R \cup \{+\infty\}$ is
a closed convex function, subdifferentiable on the closed convex
domain $\domain$.

While stochastic gradient methods are the de facto choice for
problem~\eqref{eqn:objective}---enjoying several convergence
guarantees~\cite{Zinkevich03, NemirovskiJuLaSh09, BottouBo07,
  ShalevSiSrCo11} with straightforward parallel extensions that make them
practically attractive~\cite{Lan12, DekelGiShXi12, DuchiBaWa12}---they are
sensitive to the objective $f$, noise, and hyperparameter
tuning~\cite{LiJaDeRoTa17,AsiDu19,AsiDu19siopt}.  They may even diverge
for objectives that do not satisfy their convergence criteria or with
slightly mis-specified stepsizes~\cite{AsiDu19siopt, NemirovskiJuLaSh09}.
Motivated by these limitations of gradient methods,
researchers~\cite{Bertsekas11, KulisBa10, DavisDr19, DuchiRu18c,
  AsiDu19siopt} have developed stochastic (approximate) proximal-point
(\aprox) and model-based methods as a more robust alternative.  These \aprox
methods, as we explain in Section~\ref{sec:pre}, construct a model of the
function and iterate by minimizing regularized versions of
the model. They improve over standard stochastic gradient methods, as they
are robust to stepsize choice, adaptive to problem difficulty, and converge
on a broader range of problems than stochastic gradient
methods~\cite{DuchiRu18c, AsiDu19siopt}. Yet these \aprox methods are
inherently sequential, and as we hit physical limits on processor speeds, it
is becoming clear that opportunities for improvements in large-scale
computation and energy use must focus on parallelization~\cite{FullerMi11};
it is not immediately apparent how to efficiently parallelize stochastic
model-based methods.

We study methods to parallelize the \aprox family via minibatched samples
$\statrv^{1:m} \in \statdomain^m$, that is, where each iteration of the
method receives an independent batch $\statrv^{1:m} \simiid P$, developing
several new results for model-based methods the
problem~\eqref{eqn:objective} more generally along the way. We provide the
following:
\begin{enumerate}[leftmargin=*]
\item \textit{Non-asymptotic rates and accelerated convergence:} In
  Section~\ref{sec:non-asym}, we develop nonasymptotic convergence
  guarantees that depend on the \emph{variance} of sample gradient
  estimates, in distinction to previous analyses~\cite{DavisDr19,
    AsiDu19siopt} that depend only on their magnitude, showing that
  model-based methods enjoy linear speedup in minibatch size $m$ (from
  standard $1/\sqrt{k}$ convergence rates to $1 / \sqrt{km}$), analogous to standard speedup guarantees for gradient
  methods~\cite{Lan12,DekelGiShXi12}. These also allow us to develop an
  order-optimal accelerated method for the \aprox family in
  Section~\ref{sec:acceleration}.
%% \item \textit{Acceleration:} In \Cref{sec:acceleration}, we introduce an
%%   accelerated version of the $\aprox$ method. We prove a convergence rate
%%   which is faster than the non-accelerated version, and show that this
%%   method also has speedup which scales favorably with minibatch size.
\item \label{item:interpolation} \textit{Optimal convergence and
  interpolation problems:} In Sections~\ref{sec:interpolation}
  and~\ref{sec:lower-bounds}, we consider interpolation problems, that is,
  problems for which there exists $x\opt \in \xdomain$ minimizing $F(\cdot;
  \statval)$ with $P$-probability 1. Such problems arise in numerous modern
  machine learning applications~\cite{BelkinHsMi18, BelkinRaTs19}---where
  one can achieve zero training error---or, for example, in finding a point
  in the intersection of convex sets $\cap_{i = 1}^N C_i$, where one takes
  $\statdomain = \{1, \ldots, N\}$ and $F(x; i) = \dist(x, C_i)$. For these
  problems, we both develop new optimality results, characterizing
  (worst-case) problem difficulty based on a particular growth condition Asi
  and Duchi~\cite{AsiDu19siopt} introduce, which is (by these results)
  evidently fundamental; we also give some sufficient
  conditions for minibatching to yield improved convergence.
\item
  \textit{Experimental evaluation:} We conclude with an experimental
  evaluation in Section~\ref{sec:experiments}, where we study the
  robustness and acceleration properties of the methods; performance
  profiles highlight the benefits of using these better models.
  %% In \Cref{sec:experiments}, we support
  %% our theoretical results with extensive simulation. Our experiments
  %% demonstrate that our methods with more-accurate models (e.g.,
  %% truncated~\eqref{eqn:trunc-model}) are more robust---with respect to the
  %% stepsize choice---to achieve the desirable speedups than standard
  %% stochastic gradient methods. While the latter requires a careful stepsize
  %% choice to guarantee speedup, our methods exhibit a wide range of stepsizes
  %% with significant speedup guarantees.  Moreover, as suggested by our
  %% theory, our methods exhibit an improved speedups (compared to SGM) for the
  %% non-smooth absolute regression problem. \TODO{add phase retrieval?}
\end{enumerate}

\subsection{Preliminaries}
\label{sec:pre}
The starting point of our methods is the model-based approximate
proximal-point (\aprox) framework~\cite{DavisDr19, DuchiRu18c,
  AsiDu19siopt}, which approximates the functions
$F$ via \emph{models} $F_x$ of $F$ localized at $x$,
which satisfy the following conditions:
\begin{enumerate}[label=(C.\roman*),leftmargin=*]
\item \label{cond:convex-model} \textit{Convexity:}
  The function $y \mapsto F_x(y; \statval)$ is convex
  and subdifferentiable on $\xdomain$.
\item \label{cond:lower-model} \textit{Lower bounds and local accuracy:}
  For all $y \in \xdomain$,
  \begin{equation*}
    F_x(y; \statval) \le F(y; \statval)
    ~~ \mbox{and} ~~
    F_x(x; \statval) = F(x; \statval).
  \end{equation*}
\end{enumerate}
Note that Condition~\ref{cond:lower-model} immediately implies that
$\partial F_x(y; \statval)|_{y = x} \subset \partial F(x; \statval)$.

With such a model, \aprox algorithms iteratively sample $\statrv_k \simiid
P$ and update
\begin{equation}
  \label{eqn:model-iteration}
  x_{k+1} \defeq \argmin_{x \in \xdomain}
  \left\{ F_{x_k}(x ; \statrv_k) + \frac{1}{2 \stepsize_k}
  \ltwo{x - x_k}^2 \right\}.
\end{equation}
Typical choices for the models include the following three:
\begin{itemize}[leftmargin=*]
\item \textit{Stochastic gradient methods:} for some $F'(x; \statval)
  \in \partial F(x; \statval)$, use
  the linear model
  \begin{equation}
    \label{eqn:dumb-linear-model}
    F_{x}(y; \statval) \defeq F(x; \statval) + \<F'(x; \statval), y - x\>.
  \end{equation}
\item \textit{Stochastic proximal point methods:} use the full proximal model
  \begin{equation}\label{eqn:full-prox}
    F_x(y;s) \defeq F(y; s).
  \end{equation}
\item \textit{Truncated methods:} for some
  $F'(x; \statval) \in \partial F(x; \statval)$,
  use
  \begin{equation}
    \label{eqn:trunc-model}
    F_x(y; \statval) \defeq \max\left\{
    F(x; \statval) + \<F'(x;\statval), y - x\>,
    \inf_{z \in \xdomain} F(z; \statval)\right\}.
  \end{equation}
  The model~\eqref{eqn:trunc-model} is often simple to apply: in many
  applications, the objective is non-negative, so $\inf_{z \in \xdomain}
  F(z; \statval) = 0$ and the model is simply the positive part of the
  linear approximation~\eqref{eqn:dumb-linear-model}.
\end{itemize}
%of the \aprox framework:
  
% \begin{enumerate}[label=(C.\roman*),leftmargin=*]
%   \setcounter{enumi}{2}
% \item \label{cond:lower-by-optimal}
%   For all $\statval \in \statdomain$, the models $f_x(\cdot; \statval)$ satisfy
%   %\begin{equation*}
%   $ f_x(y; \statval) \ge \inf_{z \in \xdomain}
%     f(z; \statval).$
%   %\end{equation*}
% \end{enumerate}

\paragraph{Notation}
For a convex function $f$, 
$\partial f(x)$ denotes its subgradient set at
$x$, and $f'(x) \in \partial f(x)$ denotes an arbitrary element of
the subdifferential.  
We let  $\xdomain\opt = \argmin_{x \in \xdomain} f(x)$
denote the optimal set of problem~\eqref{eqn:objective} 
and $x\opt \in \xdomain\opt$ denote a single minimizer.  
We let $\mc{F}_k \defeq \sigma(\statrv_1, \ldots,
\statrv_k)$ be the $\sigma$-field generated by the first $k$ random
variables $\statrv_i$, so $x_k \in
\mc{F}_{k-1}$ for all $k$ under iteration~\eqref{eqn:model-iteration}.

\subsection{Related work}

Stochastic gradient methods~\cite{RobbinsMo51} are the most widely used
method for solving stochastic minimization problems; an enormous literature
gives numerous convergence results~\cite{Polyak87, PolyakJu92,
  Zinkevich03, NemirovskiJuLaSh09, Zhang04, KushnerYi03, BachMo11}. The
growth of parallel computing has motivated the development of ``minibatch''
methods that use multiple samples $\statrv$ in each iteration, where
researchers have shown how stochastic gradient-like methods enjoy linear
speedups as batch sizes increase~\cite{Lan12, DekelGiShXi12, DuchiBaWa12,
  RechtReWrNi11, ChaturapruekDuRe15}. Other work proposes accelerated
stochastic optimization methods, showing faster (worst-case optimal)
associated convergence rates \cite{LinMaHa18, Lan12}.  In spite of their
successes, stochastic gradient methods still suffer a number of
drawbacks. For example, they are sensitive to problem parameters, where
mis-specified stepsizes may force slow (even order sub-optimal or
exponentially slower) convergence~\cite{NemirovskiJuLaSh09}; objective
functions without appropriate scaling or that grow too quickly may cause
divergence~\cite{AsiDu19, AsiDu19siopt}; they can fail to adapt to problem
geometry~\cite{DuchiHaSi11, LevyDu19}.  This motivates work to make
stochastic gradient methods more robust~\cite{NemirovskiJuLaSh09} and
adaptive~\cite{DuchiHaSi11, OrabonaPa18} as well as research on stochastic
proximal-point and model-based methods~\cite{KulisBa10, KarampatziakisLa11,
  Bertsekas11, DuchiRu18c, DavisDr19}.  In this vein,
Asi and Duchi~\cite{AsiDu19siopt, AsiDu19} show how better models in stochastic
optimization yield improved stability, robustness, and convergence
guarantees over classical stochastic subgradient methods.

%% Other papers recognize the
%% instability of stochastic gradients methods, demonstrating situations where
%% they can have slow convergence as a result of mis-specified
%% stepsizes~\cite{NemirovskiJuLaSh09, BachMo11, AsiDu19}, thereby confirming
%% the importance of robustness.

%% The full proximal~\eqref{eqn:full-prox} and
%% truncated~\eqref{eqn:trunc-model} models provide more accurate
%% approximations of $F$ than the linear model~\eqref{eqn:dumb-linear-model},
%% which \cite{AsiDu19siopt,AsiDu19}, motivated by previous work in the
%% area~\cite{KulisBa10, KarampatziakisLa11, DavisDr19}, show yields more
%% robust algorithms with better theoretical and practical convergence.

%% Rockafellar \cite{Rockafellar76} introduces proximal point methods, which
%% have seen a resurgence in applications to stochastic
%% optimization~\cite{KulisBa10, Bertsekas11, KarampatziakisLa11,Bianchi16,
%%   LinMaHa18}. Of most relevance to our work are extensions of the stochastic
%% proximal methods that use approximate models in the proximal
%% update~\cite{DuchiRu18c,DavisDr19,AsiDu19siopt,AsiDu19}.  Asi and Duchi
%% \cite{AsiDu19siopt} develop a stochastic approximate proximal point method,
%% namely \aprox, and establish several convergence guarantees and stability
%% properties that are superior to standard stochastic gradient methods.  Yet
%% this work does not address the challenge of minibatching and acceleration,
%% two aspects which are prevalent optimization and machine learning.
%% \TODO{Add minimax lower bounds related work? @John: We don't know the
%%   literature for this.}

A second line of work studies acceleration, mini-batching, and parallelism
in stochastic optimization~\cite{DekelGiShXi12, Lan12,
  DefazioBaLa14, NiuReReWr11, ChaturapruekDuRe15, ScamanBaBuLeMa17}. Here,
the key insights typically show that mini-batching---averaging stochastic
gradients---yields reduced variance and hence improved
convergence~\cite{DekelGiShXi12, Lan12}. Other key insights show how in
large-scale communication-limited problems, the noise inherent to sampling
dominates deterministic components of convergence rates and errors due to
delay or communication~\cite{ChaturapruekDuRe15, ManiaPaPaReRaJo17,
  ScamanBaBuLeMa17}. For model-based methods, appropriate notions of
variance are less immediate, and in interpolation problems (recall
item~\ref{item:interpolation} above) there is essentially no noise, so that
an important part of our development is to extend accelerated and
variance-dependent rates of convergence (as available for gradient-based
methods~\cite{Lan12}) to model-based methods. An important component of
accelerated and parallel methods is their (worst-case)
optimality~\cite{Nesterov04, Lan12, AgarwalBaRaWa12}; as one of the major
successes for model-based methods is in interpolation problems, it is also
of interest to develop corresponding optimality results, which (to our
knowledge) do not exist.

%% stochastic noise governs convergence rates and are

%Proximal methods with minibatching have been studied in
%\cite{PatrascuPaIr20,HendrikxBaMa19,KonecnyLiRiTa20,FercoqRi15}. They
%decompose the function into a smooth and non-smooth part, and apply
%gradient updates \wrt the smooth part and proximal updates \wrt the
%non-smooth part.

% -*- Mode: latex -*- %

\section{Methods}
\label{sec:methods}

While at some level, the extension of standard stochastic gradient methods
to parallel settings---average gradients to reduce noise---is clear, such
extension is less immediate for proximal and model-based methods.  To that
end, we identify several different possibilities for extending the \aprox
framework---which coincide for linear models (stochastic gradient
methods)---but can exhibit different optimization behavior.
Given a batch $\statrv_k^{1:\mb} \in \statdomain^{\mb}$ of samples,
we consider the following:

\textbf{Iterate averaging} (\emph{\PIA}): The naive extension of
\aprox to use minibatches is to perform
an individual update for each sample $\statrv_k^i$, then average the updates:
\begin{equation}\label{eq:iterate-avg}
  x_{k+1} \defeq \frac{1}{\mb} \sum_{i=1}^{\mb} x_{k+1}^i \quad \text{where}
  \quad x_{k+1}^i = \argmin_{x \in \mc{X}}
  \left\{ F_{x_k}(x ; \statrv_k^i) + \frac{1}{2 \stepsize_k}
  \ltwo{x - x_k}^2 \right\}.
\end{equation}
% then average these iterates to get our final update
% \begin{equation}
% \label{eq:iterate-avg}
%   .
% \end{equation}
This method's simplicity and (near) full parallelization makes it
attractive, and when $\xdomain = \R^n$ and each of the models $F_x$ is the
subgradient model~\eqref{eqn:dumb-linear-model}, it coincides with the
mini-batch stochastic gradient method.  Unfortunately, in general it does
not enjoy the same acceleration properties of our other methods.

A method that more naturally dovetails with the model-based perspective
is to minimize a model of the average
\begin{equation}
\label{eq:avg}
  \wb{F}(x ; \statrv_k^{1:m})
  \defeq \frac{1}{\mb} \sum_{i=1}^{\mb} 
     F(x ;\statrv_k^i)
\end{equation}
at every iteration. In particular, with any model
$\Fbar_{x_k}(x;S_k^{1:\mb})$ of the average satisfying Conditions
\ref{cond:convex-model} and \ref{cond:lower-model}, we can perform the
update
\begin{equation}
  \label{eq:model-avg}
  x_{k+1} \defeq \argmin_{x \in \mc{X}}
  \left\{ \wb{F}_{x_k}(x ; \statrv_k^{1:m}) 
  + \frac{1}{2 \stepsize_k} \ltwo{x - x_k}^2 \right\}.
\end{equation}
While our theorems hold for any model-based algorithm satisfying Conditions
\ref{cond:convex-model} \ref{cond:lower-model} (and
Condition~\ref{cond:trunc} to come), we find two instatiations of the
approach~\eqref{eq:model-avg} of particular interest.

\newcommand{\lowerbound}{\Lambda}

\textbf{Truncated Average} (\emph{\PMA}): The first such model
extends the truncated
model~\eqref{eqn:trunc-model}. Let $\lowerbound(\statval^{1:\mb})$ be any
lower bound on $\wb{F}(\cdot, \statval^{1:\mb})$; for example,
$\lowerbound(\statval^{1:\mb}) = \frac{1}{\mb} \sum_{i=1}^\mb \inf_{z \in
  \mc{X}} F(z; \statval^i)$ suffices.  Then set
\begin{equation*}
  \wb{F}_x(y; \statval^{1:m})
  \defeq \max\Big\{\wb{F}(x; \statval^{1:\mb})
  + \<\wb{F}'(x; \statval^{1:\mb}), y - x\>,
  \lowerbound(\statval^{1:\mb})
  \Big\}.
\end{equation*}
In the standard case that the functions $F$ are nonnegative and $\xdomain =
\R^n$, the update~\eqref{eq:model-avg} corresponds to (stochastic) Polyak
stepping~\cite{Polyak87}, and becomes
\begin{equation}
  \label{eq:model-avg-trunc}
  x_{k+1} = x_k - \min \left\{\stepsize_k, 
  \frac{\wb{F}(x_k; \statrv_k^{1:m}) - \lowerbound(S^{1:m})}{\ltwos{\wb{F}'(x_k; \statrv_k^{1:m})}^2} \right\} \wb{F}'(x_k; \statrv_k^{1:m}).
\end{equation}
% \TODO{shouldn't this update depend on what we choose to be the lower bound.}
The update~\eqref{eq:model-avg-trunc} for the truncated models thus yields
an embarrassingly parallelizable scheme: each worker computes $F(x_k
;\statrv_k^i)$ and $ \nabla F(x_k ;\statrv_k^i)$, which need only be
averaged to apply the update~\eqref{eq:model-avg-trunc}. 
% As we show later, this approach provides several robustness and convergence guarantees.

\textbf{Average of Truncated Models} (\emph{\PAM}): The
update~\eqref{eq:model-avg-trunc} ignores some structural aspects of the
objectives $F$; it is natural to consider a more accurate averaging
of models. Letting
$F_{x}(y;\statval^i) = \max\{F(x; \statval^i) + \<F'(x; \statval^i), y -
x\>, \inf_{z \in \domain} F(z;\statval^i)\}$, the average
$\frac{1}{\mb} \sum_{i = 1}^{\mb} F_{x}(\cdot; \statval^i)$ satisfies
conditions~\ref{cond:convex-model} and~\ref{cond:lower-model}, and we
consider the update
\begin{equation}
\label{eq:avg-model}
  x_{k+1} \defeq \argmin_{x \in \mc{X}}
  \left\{\frac{1}{\mb} \sum_{i=1}^\mb F_{x_k}(x;\statrv_k^i) 
  + \frac{1}{2 \stepsize_k} \ltwo{x - x_k}^2 \right\}.
\end{equation}
When $\mb$ is not too large, problem~\eqref{eq:avg-model}
is relatively easy to solve. Indeed, define $g_i = F'(x_k; \statrv_k^i)$ and
let $G = [g_1 ~ \cdots ~ g_\mb] \in \R^{n \times m}$ and
$v = [F(x_k; \statrv_k^1) ~ \cdots ~ F(x_k; \statrv_k^\mb)]^T \in \R^\mb$.
Then the dual to problem~\eqref{eq:avg-model} is
\begin{equation*}
  \begin{split}
    \maximize ~ & -\frac{\stepsize}{2} \lambda^T G^T G \lambda + \lambda^T v
    ~~~~
    \subjectto ~ 0 \preceq \lambda \preceq \frac{1}{m},
  \end{split}
\end{equation*}
and letting $\lambda_k$ be the solution, we update $x_{k + 1} = x_k -
\stepsize_k G \lambda_k$.  In situations where computing the
(sub)gradients $F'(x_k; \statrv_k^i)$ or infima $\inf_z F(z; \statrv_k^i)$
are more expensive than solving the dual---reasonable when $m$ is
small---one can parallelize efficiently.

\begin{remark} The preceding two methods provide two natural approaches
  to mini-batching model-based stochastic methods; any approach that
  guarantees the model $\wb{F}_x(y; \statval^{1:\mb})$ satisfies
  conditions~\ref{cond:convex-model}--\ref{cond:lower-model} for
  the average $\wb{F}(\cdot; \statval^{1:\mb})$ will similarly suffice
  for our results.
  %% It We note that we could easily generalize \PMA~and
  %% \PAM~methods beyond the truncated model. Let $\mc{H}$ denote the set of
  %% convex functions. Let $G: \mc{H} \times \domain \rightarrow \mc{H}$ be some
  %% operator which takes in a convex function $F$ and a point $y$ in the domain
  %% as input and constructs a model of $F$ at $y$ satisfying Conditions
  %% \ref{cond:convex-model} and \ref{cond:lower-model}. Then the \PMA~model is
  %% more generally written as $G(\Fbar(x; S^{1:m}), x_k)$. The \PAM~model is
  %% more generally $\sum_{i=1}^m G(F(x; S^i), x_k)/m$.
\end{remark}

\begin{figure}[ht]
  \begin{center}
    \begin{tabular}{ccc}
      \begin{overpic}[width=.3\columnwidth]{%,grid]{%
      		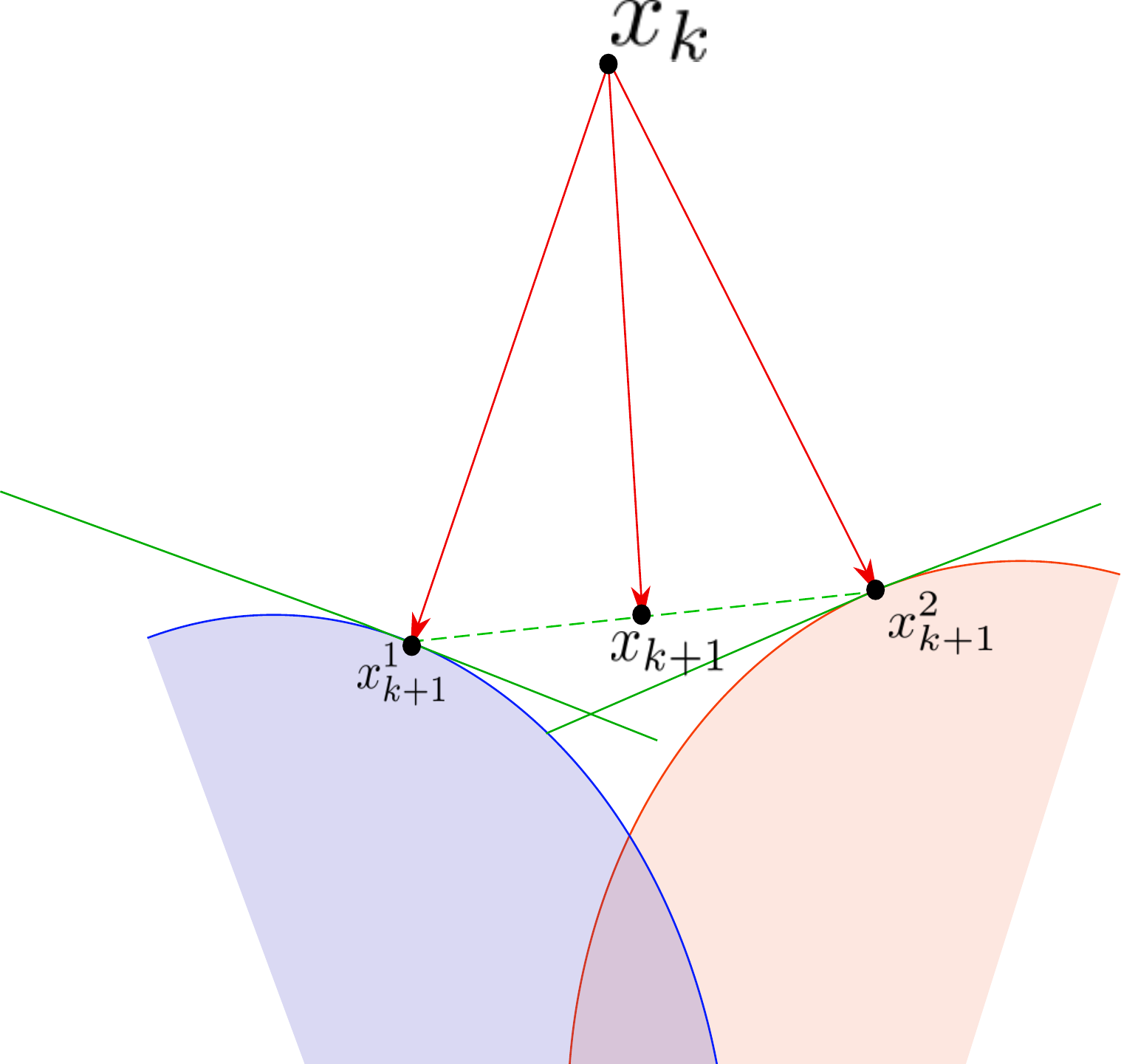}
      \end{overpic} &
      \begin{overpic}[width=.3\columnwidth]{%,grid]{%
      		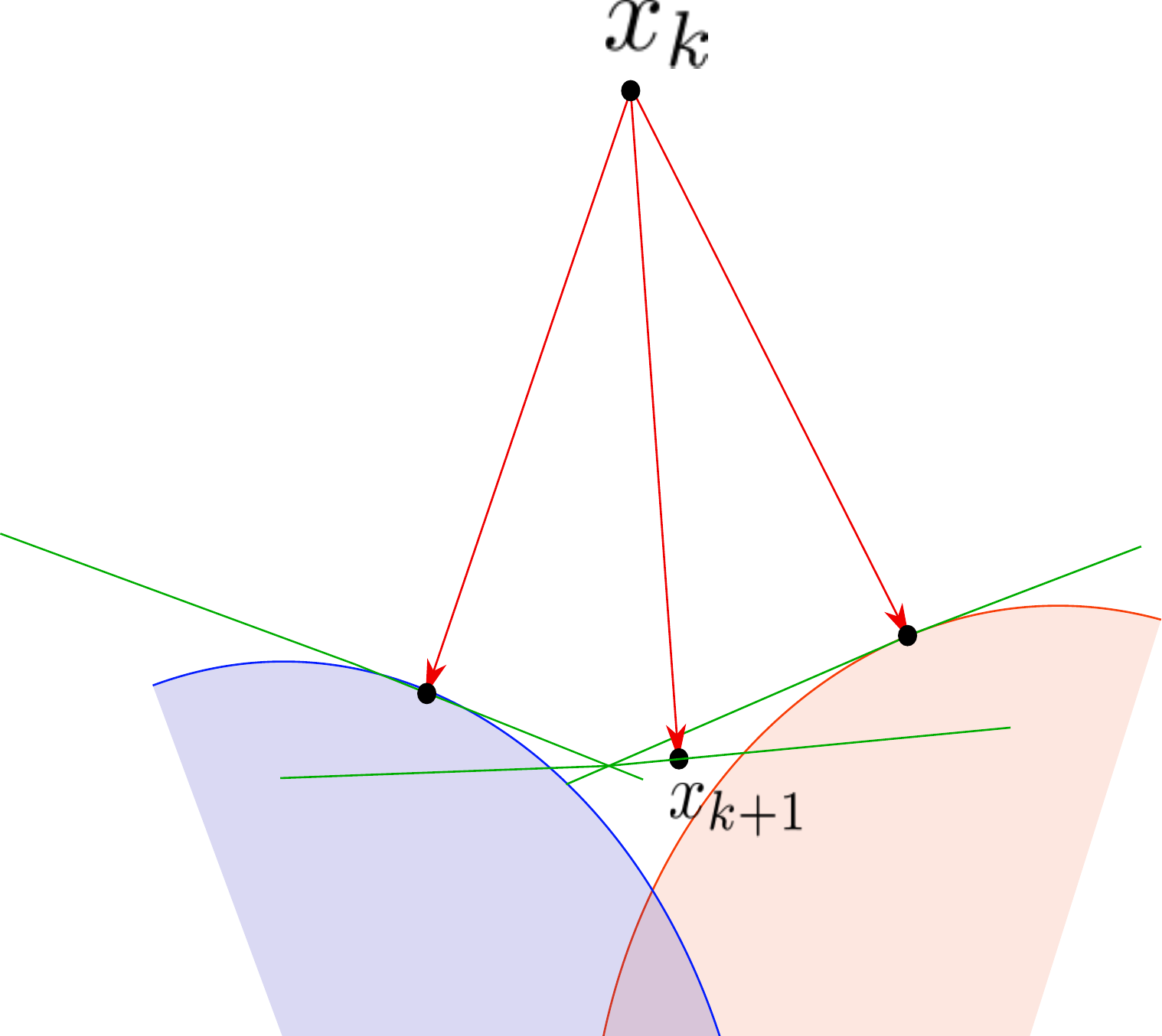}
      \end{overpic} &
      \begin{overpic}[width=.3\columnwidth]{%,grid]{%
      		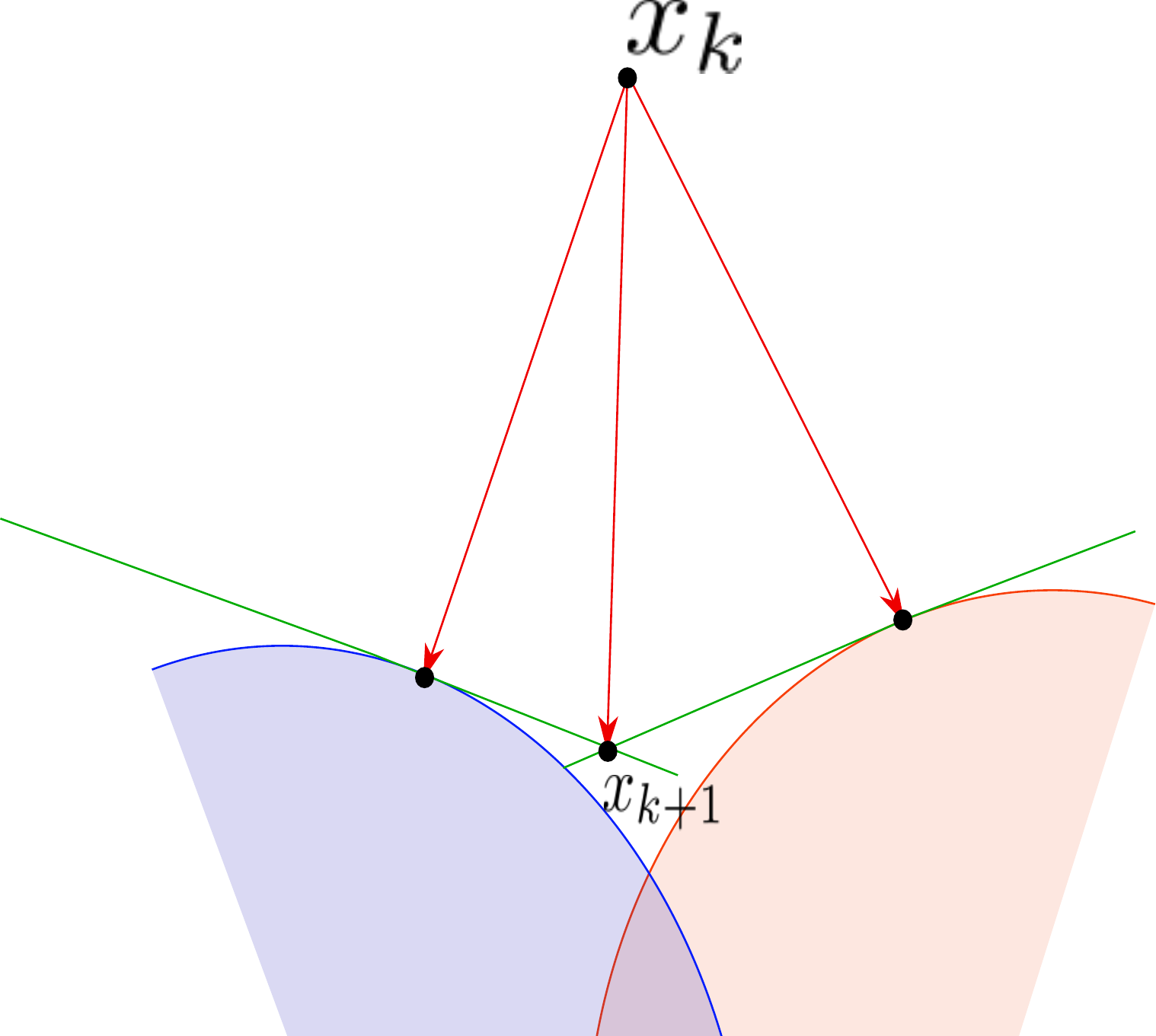}
      \end{overpic} \\
  	  (a) & (b) & (c)
    \end{tabular}
  \caption{
    \label{fig:methods-intuition}
    Updates for the truncated model
    using (a) iterate averaging~\eqref{eq:iterate-avg},
    (b) truncated averaging~\eqref{eq:model-avg-trunc},
    and (c) averaging of models~\eqref{eq:avg-model}.
  }
  \end{center}
  \vspace{-.5cm}
\end{figure}

Before proceeding to our theoretical guarantees, we provide a simple example
to help illustrate the methods. Consider the problem of finding a point in
the intersection $C_1 \cap C_2$ of convex sets $C_1$ and $C_2$ by minimizing
$f(x) = \half(\dist(x, C_1) + \dist(x, C_2))$.
Figure~\ref{fig:methods-intuition} illustrates the \PIA, \PMA, and
\PAM~updates given infinite stepsize $\stepsize$ (which still guarantees
convergence if $\cap_i C_i$ is non-empty~\cite{AsiDu19siopt}). Let
$\pi_i(x)$ denote the projection of $x$ onto $C_i$, so that the hyperplane
tangent to $C_i$ at $\pi_i(x)$ is $x - \pi_i(x)$.  In this case, iterate
averaging~\eqref{eq:iterate-avg} projects the current iterate $x_k$ to the
two sets in the batch and averages these updates
(Fig.~\ref{fig:methods-intuition}(a)).  The \PMA
update~\eqref{eq:model-avg-trunc} constructs the average of the hyperplanes
$v = \half \sum_{i = 1}^2 (x - \pi_i(x))$ and distances $d =
\half \sum_{i = 1}^2 \dist(x, C_i)$ and projects to
the halfspace $\{x \mid v^T (x - x_k) \le -d\}$, which yields more progress
(Fig.~\ref{fig:methods-intuition}(b)).
The \PAM update~\eqref{eq:avg-model} projects to the set defined by the
intersection of the tangent halfspaces at $\pi_i(x)$
(Fig.~\ref{fig:methods-intuition}(c)). We expect generally---and
our theory and experiments will confirm---that the
\PMA and \PAM updates are more effective than naive averaging.

\section{Non-Asymptotic Convergence Results} \label{sec:non-asym}

%In this section, we discuss the non-asymptotic convergence guarantees of {\it model of average} \eqref{eq:model-avg} and {\it average of models} \eqref{eq:avg-model} methods as described in the \Cref{sec:methods}. In particular, we show that, under smoothness, we get an $\mb$-factor speed up in the convergence rate of these methods. The idea is to leverage the fact that variance of the gradients reduces by a factor of $\mb$ on averaging. For this, we first prove a result for generic model based updates where the model satisfies conditions C1 and C2 {\color{red} add proper ref to C1 and C2}. This highlights the variance of gradient in the computational complexity. Then, by showing the reduction in variance due to averaging we show the speedup obtained due to minibatching.

\providecommand{\lipgrad}{L}
\newcommand{\divergence}{D_\distgen}
\newcommand{\distgen}{h}  % Distance generating function
\newcommand{\fnerror}{e}
\newcommand{\fnerrory}{e^{(1)}}
\newcommand{\graderror}{\xi}
\newcommand{\radius}{R}
\newcommand{\gradvariance}{\sigma_0}  % Gradient variance
\newcommand{\linapprox}{\mathsf{lin}_f}  % Linear approximation

Our first set of theoretical results extends the familiar
non-asymptotic rates of convergence for smooth convex stochastic
optimization~\cite{Lan12} to model-based methods. Here, we show that
model-based methods for problem~\eqref{eqn:objective} enjoy optimal
dependence on the variance of stochastic gradients, and, building off of
Tseng~\cite{Tseng08} and Lan~\cite{Lan12}, can be accelerated to achieve
worst-case optimal complexity.  To present our results in the most
generality, we allow non-Euclidean geometries to generalize mirror
descent~\cite{BeckTe03, NemirovskiJuLaSh09}.

To that end, recall that a differentiable convex function $\distgen$ is a
\emph{distance generating function} for $\xdomain$ if it is strongly convex
with respect to a norm $\norm{\cdot}$ over $\xdomain$, meaning $\distgen(y)
\ge \distgen(x) + \<\nabla \distgen(x), y - x\> + \half \norm{x - y}^2$ for
$x, y \in \xdomain$. The associated \emph{Bregman divergence} is then
$\divergence(x, y) \defeq h(x) - h(y) - \<\nabla h(y), x - y\>$, which
evidently satisfies $\divergence(x, y) \ge \half \norm{x - y}^2$.  Recalling
the dual norm $\dnorm{z} = \sup_{\norm{x} \le 1} \<z, x\>$, throughout this
section, we will work with the following standard assumption~\cite{Lan12,
  DekelGiShXi12}.
\begin{assumption}
  \label{assumption:bounded-noise}
  The function $f$ has $\lipgrad$-Lipschitz gradient with respect
  to the norm $\norm{\cdot}$, meaning that
  \begin{equation*}
    \dnorm{\nabla f(x) - \nabla f(y)} \le \lipgrad \norm{x - y},
  \end{equation*}
  and there exists
  $\gradvariance^2 < \infty$ such that for each $x \in \xdomain$,
  \begin{equation*}
    \E[\dnorm{\nabla f(x) - \nabla F(x;\statrv)}^2] \leq \gradvariance^2.
  \end{equation*}
\end{assumption}

When $\divergence(x, y) \le \radius^2$ for all $x, y \in \xdomain$ and
Assumption~\ref{assumption:bounded-noise} holds, mirror descent methods
achieve convergence guarantees of the form $\frac{\lipgrad \radius^2}{k} +
\frac{\gradvariance \radius}{\sqrt{k}}$, while accelerated
methods~\cite{Lan12} can achieve $\frac{\lipgrad \radius^2}{k^2} +
\frac{\gradvariance \radius}{\sqrt{k}}$. The latter is worst-case
optimal~\cite{NemirovskiYu83}.
By considering the natural generalization
\begin{equation}
  \label{eqn:mirror-model}
  x_{k + 1} = \argmin_{x \in \xdomain}
  \left\{F_{x_k}(x; \statrv_k) + \frac{1}{\stepsize_k} \divergence(x, x_k)
  \right\}
\end{equation}
of the model-based iteration~\eqref{eqn:model-iteration}, we achieve
the same
(optimal) rates here for the model-based mirror
method~\eqref{eqn:mirror-model}; combined with the results of the
paper~\cite{AsiDu19siopt}, these show that model-based methods offer the
same benefits (efficiency, parallelizability, and worst-case optimality)
that stochastic gradient and mirror descent methods do while simultaneously
guaranteeing more robustness.

%% We prove variance-based convergence results first for the non-accelerated algorithm and then for the accelerated version. Rates depending on the minibatch size follow directly from the variance-based results because minibatching reduces the variance of the gradient estimate. 
%% We begin in \Cref{sec:smooth} 
%% with a standard speedup result, resembling known results for SGM. 
%% In \Cref{sec:non-smooth} we show that our methods using the 
%% full proximal model~\eqref{eqn:full-prox} enjoy
%% linear speedup in minibatch size even for non-smooth
%% functions. Finally in \Cref{sec:acceleration}, we show convergence results for the accelerated algorithm.

\subsection{A basic non-asymptotic convergence guarantee}
\label{sec:smooth}

Our first result gives convergence of the basic
iteration~\eqref{eqn:mirror-model}.
\begin{theorem}
  \label{theorem:non-asymptotic}
  Let Assumption~\ref{assumption:bounded-noise} hold, and assume
  that $\divergence(x, y) \le \radius^2$ for all $x, y \in \xdomain$.
  Let
  $x_k$ follow the model-based
  iteration~\eqref{eqn:model-iteration} for any model
  satisfying Conditions~\ref{cond:convex-model} and~\ref{cond:lower-model}.
  Define the stepsizes $\stepsize_k = \frac{1}{\lipgrad + \eta_k}$, where
  $\eta_k \ge 0$ is non-decreasing. Then
  \begin{equation*}
    \sum_{i = 1}^{k}
    \E[f(x_{i+1}) - f(x\opt)]
    \leq \frac{\lipgrad \radius^2}{2} +
    \frac{\radius^2 \eta_k}{2} + \sum_{i = 1}^k \frac{\gradvariance^2}{2\eta_i}.
  \end{equation*}
\end{theorem}
\noindent
The proof of this result, while not completely standard as we cannot rely
on linearity in the updates or gradients to achieve the variance
bound, builds off of several well-established techniques, so we defer
it to Appendix~\ref{sec:proof-non-asymptotic}.

Having established a convergence result that
depends on the noise of the gradient estimates, convergence
guarantees for the average $\wb{x}_k = \frac{1}{k}
\sum_{i = 1}^k x_{i + 1}$ are immediate. First,
under the conditions of Theorem~\ref{theorem:non-asymptotic}
we have
\begin{equation*}
  \E[f(\wb{x}_k)] - f(x\opt) \le \frac{\lipgrad \radius^2}{k}
  + \frac{\radius^2 \eta_k}{k} + \frac{\gradvariance^2}{2 \eta_k},
\end{equation*}
and with the choice $\eta_k = \eta_0 \gradvariance \sqrt{k} / \radius$ we obtain the
rate
\begin{equation}
  \label{eqn:naive-rate}
  \E[f(\wb{x}_k)] - f(x\opt) \le \frac{\lipgrad \radius^2}{k}
  + \frac{\radius \gradvariance}{\sqrt{k}}
  \left(\eta_0 + \frac{1}{2 \eta_0}\right).
\end{equation}
When we use the standard Euclidean choice
$\distgen(x) = \half \ltwo{x}^2$, we see an immediate speedup guarantee
for the minibatched \aprox methods:
\begin{corollary}\label{cor:non-asym-variance}
  Let the conditions of \Cref{theorem:non-asymptotic} hold, let $\eta_k =
  \eta_0 \sqrt{k}$ with $\eta_0 = \frac{\gradvariance}{\sqrt{\mb} \radius}$, and
  let $x_k$ be generated by the iteration~\eqref{eq:model-avg} with any
  model $\Fbar_{x}(y;\statrv^{1:m})$ satisfying
  conditions~\ref{cond:convex-model} and~\ref{cond:lower-model} and
  minibatch size $\mb$.  Then
  \begin{equation*}
    \E[f(\wb{x}_k) - f(x\opt)]
    \leq \frac{\lipgrad \radius^2}{k} + \frac{3 \radius\gradvariance}{2\sqrt{k \mb}}.
  \end{equation*}
\end{corollary}
When the iteration count $k \gg \frac{\lipgrad^2
  \radius^2\mb}{\gradvariance^2}$, the second term dominates the rate of
convergence. Letting $T(\epsilon)$ denote the number of iterations to
achieve $\E[f(\wb{x}_{T(\epsilon)}) - f(x\opt)] \le \epsilon$, we obtain
that $T(\epsilon) \lesssim \frac{R^2 \gradvariance^2}{\epsilon^2 m}$, that is,
there is a linear speedup as a function of the minibatch of size $\mb$.
This is similar to the speedup that standard stochastic gradient methods
achieve~\cite{Lan12,DekelGiShXi12} and is minimax optimal.

\subsection{Accelerated model-based methods}
\label{sec:acceleration}

\newcommand{\lotserror}{\zeta}

We now develop an accelerated analogue of the
iteration~\eqref{eqn:model-iteration}, which gives a leading minimax-optimal
$O(1/k^2)$ rate, by building off of the ideas of Lan~\cite{Lan12} and
Tseng~\cite{Tseng08}. We consider a modified iteration, which
augments the model-based update~\eqref{eqn:model-iteration} with two
auxiliary sequences whose momentum allows accelerated convergence.
For full generality and completeness, we consider an augmented version
of problem~\eqref{eqn:objective}, where we wish to minimize
\begin{equation*}
  f(x) + r(x) = \E_P[F(x; \statrv)] + r(x),
\end{equation*}
where $r$ is a known convex function (typically a regularizer of some type).
We require a  non-increasing sequence $\theta_k \in [0, 1]$ of
stepsizes and consider the three term iteration
\begin{equation}
  \label{eqn:accelerated}
  \begin{split}
    y_k & = (1 - \theta_k) x_k + \theta_k z_k \\
    z_{k + 1} & = \argmin_{x \in \xdomain}
    \left\{ F_{y_k}(x; \statrv_k) + r(x) + \frac{1}{\stepsize_k}
    \divergence(x, z_k) \right\} \\
    x_{k + 1} & = (1 - \theta_k) x_k + \theta_k z_{k + 1}.
  \end{split}
\end{equation}
All our analysis requires is that the additional stepsizes $\theta_k$
satisfy $\theta_0 = 1$,
%% \begin{equation*}
$\frac{1 - \theta_k}{\theta_k^2} \le \frac{1}{\theta_{k-1}^2}$
%% \end{equation*}
for all $k$, and are non-increasing; for example, our
choice $\theta_k =
\frac{2}{k + 2}$ satisfies these desiderata, as does taking any
constant stepsize.
We then have the following theorem.
\begin{theorem}
  \label{theorem:accelerated}
  Let Assumption~\ref{assumption:bounded-noise} hold, and assume that
  $\divergence(x\opt, x) \le \radius^2$ for all $x \in \xdomain$.  Let $(y_k,
  z_k, x_k)$ follow the three term iteration~\eqref{eqn:accelerated} for any
  model satisfying Conditions~\ref{cond:convex-model}
  and~\ref{cond:lower-model}. Take stepsizes
  $\theta_k = \frac{2}{k + 2}$ and
  $\stepsize_k = \frac{1}{\lipgrad + \eta_k}$ for
  $\eta_k = \eta_0 \sqrt{k + 1}$, where $\eta_0 \ge 0$. Then
  \begin{equation*}
    \E[f(x_{k + 1}) + r(x_{k + 1}) - f(x\opt) - r(x\opt)]
    \le \frac{4 \lipgrad \radius^2}{(k + 2)^2}
    + 2 \frac{\radius^2}{\sqrt{k}} \left[\frac{\gradvariance^2}{\eta_0} + \eta_0
      \right].
  \end{equation*}
\end{theorem}
\noindent
See Appendix~\ref{sec:proof-accelerated} for a proof.

Specializing to the ``minibatch'' setting with
$\distgen(x) = \half \ltwo{x}^2$ again yields a minimax optimal
algorithm for the class of problems we consider.
\begin{corollary}
  Let the conditions of Theorem~\ref{theorem:accelerated}
  hold, except that we use a minibatch $\statrv_k^{1:\mb} \simiid P$ of
  size $\mb$ at each iteration, and $\wb{F}_{y_k}(\cdot; \statrv_k^{1:\mb})$
  is a model of $\frac{1}{m} \sum_{i = 1}^m F(\cdot; \statrv_k^i)$ satisfying
  Conditions~\ref{cond:convex-model} and~\ref{cond:lower-model}.
  Set $\eta_0 = \frac{\gradvariance \sqrt{\mb}}{\radius}$. Then
  \begin{equation*}
    \E[f(x_{k + 1}) + r(x_{k + 1}) - f(x\opt) - r(x\opt)]
    \le \frac{4 \lipgrad \radius^2}{(k + 2)^2}
    + 3 \frac{\radius \gradvariance}{\sqrt{k \mb}}.
  \end{equation*}
\end{corollary}

\noindent
The error rate $\bigO{1/k^2 + 1/\sqrt{k \mb}}$ is faster than the $\bigO{1/k
  + 1/\sqrt{k\mb}}$ rate we showed for the basic minibatched \aprox
algorithm~\eqref{eqn:model-iteration}, and it is minimax rate optimal.

\section{Interpolation Problems}
\label{sec:interpolation}

In \emph{interpolation problems}, there exists a consistent solution $x\opt
\in \xdomain$ satisfying $F(x\opt; \statrv) = \inf_{z \in \xdomain} F(z;
\statrv)$ with probability 1. While this is a strong assumption, it holds in
numerous practical scenarios: in machine learning problems, where a perfect
predictor (at least on training data) exists~\cite{BelkinHsMi18,
  BelkinRaTs19, MaBaBe18}; in problems of finding a point in the
intersection $C\opt = \cap_{i = 1}^N C_i$ of convex sets $C_i$, assuming
$C\opt \neq \emptyset$, where we may take $F(x; i) = \dist(x, C_i)$
(e.g.~\cite{BauschkeBo96}); or in least-squares problems with consistent
solutions~\cite{NeedellWaSr14, StrohmerVe09}. We show
a few results in this section, first that model-based methods
(often) enjoy linear convergence on these problems---in analogy
to the results available for stochastic gradient
methods~\cite{MaBaBe18}---while also demonstrating improvement via
mini-batching and reducing variance. Second, we revisit
the convergence guarantees that Asi and Duchi~\cite{AsiDu19siopt}
provide, giving a unified treatment and some discussion of the possibilities
of parallelism. These conditions appear on their face to be somewhat
non-standard, but as we show, they capture the essential difficulty
of interpolation problems, and we can provide sharp (matching to within
numerical constants) lower bounds for optimization using them.

%% We consider the following class of problems.
\begin{definition}
  \label{def:interpolation}
  Let $\xdomain\opt \defeq \argmin_{x \in \xdomain} f(x)$. Then
  problem~\eqref{eqn:objective} is an \emph{interpolation problem} if there
  exists $x\opt \in \xdomain\opt$ such that for $P$-almost all $\statval \in
  \statdomain$, we have $\inf_{x \in \domain} F(x;\statval) = F(x\opt ;
  \statval)$.
\end{definition}

%% While this restricts the class of objectives, many problems satisfy
%% this condition~\cite{MaBaBe18, SchmidtLe13, AsiDu19siopt}, including
%% overdetermined linear systems, finding a point in the intersection of convex
%% sets, and the modern machine learning problems where it is possible to
%% achieve zero training loss~\cite{BelkinHsMi18, BelkinRaTs19}.
%For example consider solving an overdetermined linear system (i.e., $\min_{x\in \domain}$ \TODO) or finding a point in the intersection of collection of convex sets (i.e., $\min{x \in \domain}$ \TODO). See \TODO for a more detailed elaboration on examples of such problems. 

We develop two sets of upper bounds for such interpolation problems.  The
first applies to any model-based method, while the second relies on the
models having more fidelity to the functions $F$.

\subsection{Upper bounds under smoothness and quadratic growth}
\label{sec:upper-bounds}

Our first set of upper bounds relies on two assumptions about the growth
of the function $f$ at the optimum---which is weaker than
typical strong convexity assumptions~\cite{MaBaBe18} that require quadratic
growth everywhere---and the noise in its gradients.
\begin{assumption}[Quadratic Population Growth]
  \label{assm:gamma-growth-population}
  There exist $\lambda> 0$ such that for all $x\in \domain$,
  \begin{align*}
    f(x) - f(x\opt)\geq \lambda \dist(x, \domain\opt)^2.
  \end{align*}
\end{assumption}
\begin{assumption}
  \label{assm:bounded-noise-dist}
  There exists $\vargrowth^2 < \infty$ such that for every $x \in \domain$,
  we have $\E[\ltwo{f'(x) - F'(x;\statrv)}^2] \leq \vargrowth^2 \dist(x,
  \domain\opt)^2$.
\end{assumption}

\noindent
It is straightforward to give examples satisfying the assumptions; noiseless
linear regression problems provide the simplest such approach.

\begin{example}
  Consider a linear regression problem with data $\statval = (a, b) \in \R^n
  \times \R$, where $a^T x\opt = b$ for all $(a, b)$, and set $F(x; (a,b)) =
  \half (a^T x - b)^2$.
  If the data $a$ belong to a subspace $V \subset \R^n$ (which may
  be $V = \R^n$), then
  Assumption~\ref{assm:gamma-growth-population} holds with
  $\lambda = \inf_{\ltwo{v} = 1} \{v^T \E[aa^T] v / 2 \mid v \in V\}$,
  and it is immediate that $\var( F'(x; \statrv)) \le
  \E[\ltwo{a}^2 \<a, x - x\opt\>^2]$, so
  Assumption~\ref{assm:bounded-noise-dist} holds with $\vargrowth^2 =
  \lambda_{\max}(\E[\ltwo{a}^2aa^T])$.
  For example, if $a$ is uniform on the scaled sphere $\sqrt{n} \sphere^{n-1}$,
  then $\lambda = 1$ and $\vargrowth^2 = n$.
\end{example}

\noindent
Alternatively, we may follow Ma et al.~\cite{MaBaBe18} by considering
a problem where the functions $F$ have Lipschitz gradients:

\begin{example}
  \label{example:smooth-instantaneous}
  If  $F(\cdot; \statval)$ has $\lipgrad(\statval)$-Lipschitz
  gradient and problem~\eqref{eqn:objective} is an interpolation problem
  with $x\opt \in \interior \xdomain$, then $\nabla F(x\opt; \statrv) = 0$
  with probability 1, and so
  \begin{align*}
    \E[\ltwo{\nabla f(x) - \nabla F(x; \statrv)}^2]
    & = \E[\ltwo{\nabla f(x) - \nabla f(x\opt)
        - (\nabla F(x\opt; \statrv) - \nabla F(x; \statrv))}^2] \\
    & \le 2 \ltwo{\nabla f(x) - \nabla f(x\opt)}^2
    + 2 \E[\ltwo{\nabla F(x\opt; \statrv) - \nabla F(x; \statrv)}^2] \\
    & \le 4 \E[\lipgrad(\statrv)^2] \ltwo{x - x\opt}^2.
  \end{align*}
  We may thus take $\vargrowth^2 \lesssim \E[\lipgrad(\statrv)^2]$.
\end{example}

Under these assumptions, model-based methods enjoy
linear (or nearly linear) convergence with constant and decaying
stepsize choices.

\begin{theorem}\label{thm:easygammagrowth}
  Assume problem~\eqref{eqn:objective} is an interpolation problem
  (Definition~\ref{def:interpolation}) and let $f$ have $\smooth$-Lipschitz
  gradient and satisfy Assumptions~\ref{assm:gamma-growth-population} and \ref{assm:bounded-noise-dist}
  , where $\lipgrad \ge \lambda$.
  Let $x_k$ follow the model-based
  iteration~\eqref{eq:model-avg} with any model $\Fbar_{x}(y; \statrv^{1:m})$
  satisfying conditions \ref{cond:convex-model} and \ref{cond:lower-model}
  with minibatch size $\mb$. Then
  \begin{enumerate}[label = (\roman*), leftmargin=0.75cm]
    \item Let $\stepsize_k = \frac{1}{\lipgrad + \eta_k}$ for $\eta_k \ge
      0$. Then
      \begin{align*} 
        \E[\dist(x_k,\domain\opt)^2]
        & \leq \exp\left(
        - \half \sum_{i = 1}^k \lambda \stepsize_k
        + \sum_{i = 1}^k \frac{\vargrowth^2}{\mb} \frac{\stepsize_i}{\eta_i}
        \right) \dist(x_0, \xdomain\opt)^2.
      \end{align*}
    \item With the constant stepsize choice
      $\stepsize_k = (\lipgrad + \eta)^{-1}$ and $\eta =
      \max\{L,\frac{8\vargrowth^2}{m\lambda}\}$,
      \begin{align*}
        \E[\dist(x_k,\domain\opt)^2] &\leq \exp\left(- k\min\left\{\frac{\scvx}{8L},\frac{m\lambda^2}{64\vargrowth^2}\right\} \right)\E[\dist(x_0,\domain\opt)^2].
      \end{align*}
  \end{enumerate}
\end{theorem}
%%%%%%%%%%%%%%%%%%%%%%%%%%%%%%%%%%%%%%%%%%%%%%%%%%%%%%
\begin{proof}
  We assume without loss of generality that $f(x\opt) = 0 = F(x\opt;
  \statval)$ for notational simplicity.  We begin with the single step
  guarantee of Lemma~\ref{lemma:smooth-progress}. Let $D_k = \dist(x_k,
  \xdomain\opt)$ for shorthand, and recall our notations $\fnerror_k =
  (F(x\opt; \statrv_k) - f(x\opt)) - (F(x_k; \statrv_k) - f(x_k)) = f(x_k) -
  F(x_k; \statrv_k)$ (in this case) and $\graderror_k = \nabla F(x_k;
  \statrv_k) - \nabla f(x_k)$. Then Lemma~\ref{lemma:smooth-progress} implies
  \begin{align*}
    \half D_{k+1}^2
    & \le \half D_k^2 - \stepsize_k f(x_{k+1})
    + \stepsize_k \fnerror_k + \frac{\stepsize_k}{2\eta_k}\ltwo{\graderror_k}^2 \\
    & \le
    \half D_k^2 - \frac{\stepsize_k \lambda}{2} D_{k+1}^2
    + \stepsize \fnerror_k +\frac{\stepsize_k}{2 \eta_k}\ltwo{\graderror_k}^2,
  \end{align*}
  where the second inequality follows from
  Assumption~\ref{assm:gamma-growth-population} that $f(x_{k+1})\geq
  \frac{\lambda}{2} D_{k+1}^2$.
  Noting that $\E[\ltwo{\graderror_k}^2 \mid x_k] \le \frac{\vargrowth^2}{\mb}
  D_k^2$ by Assumption~\ref{assm:bounded-noise-dist},
  we rearrange and take expectations on both
  sides to obtain
  \begin{align*}
    \E[D_{k+1}^2] &\leq   \underbrace{\frac{1}{\stepsize_k \lambda + 1}}_{\le \exp(-\stepsize_k \lambda/2)}\underbrace{\left(1 + \frac{\stepsize_k \vargrowth^2}{\eta \mb} \right)}_{\le \exp(\frac{\stepsize_k \vargrowth^2}{\eta \mb})}\E[D_k^2]
    \le \exp\left(\frac{-\scvx \stepsize_k}{2}
    + \frac{\vargrowth^2 \stepsize_k}{\eta_k \mb} \right) \E[D_k^2].
  \end{align*}
  Iterate this inequality to achieve the result (i) in the theorem.

  For result (ii), we simply note that if $\stepsize_k = \frac{1}{\lipgrad +
    \eta}$, then using $2\max\{L,\eta\}> L + \eta
  > \eta$, we have
  \begin{align*}
    \E[D_{k+1}^2] &\leq \exp\left(\frac{-\scvx}{4\max\{L,\eta\}} + \frac{\vargrowth^2}{\eta^2\mb} \right)
    \E[D_k^2].
  \end{align*}
  Substituting $\eta = \max\{L,\frac{8\vargrowth^2}{m\lambda}\}$ gives the
  result.
\end{proof}
% ~\\
% \noindent {\it Case 1:} $\mb > 8\vargrowth^2/(\scvx\smooth)$:\\
% Choosing $\eta = L$, we have 
% \begin{align*}
% \mathbb{E}[D_{k+1}^2] &\leq \exp\left(\frac{-\scvx}{4L} + \frac{\scvx}{8L} \right)
% \E[D_k^2] = \exp\left(\frac{-\scvx}{8L} \right)\E[D_k^2].
% \end{align*}

% \noindent {\it Case 2:} $\mb \le 8\vargrowth^2/(\scvx\smooth)$:\\
% Choosing $\eta = \frac{8\sigma_2^2}{m\lambda}$, we have 
% \begin{align*}
% \mathbb{E}[D_{k+1}^2] &\leq \exp\left(\frac{-m\scvx^2}{64\vargrowth^2}\right)\E[D_k^2].
% \end{align*}
% Reapplying these inequality repeatedly completes the proof.

The results in Theorem~\ref{thm:easygammagrowth} imply that when
the batch size is large enough that $\mb \gtrsim \vargrowth^2 / (\scvx
\smooth)$, we achieve
convergence rate $\E[\dist(x_k,\domain\opt)] \lesssim (1 - c
\frac{\lambda}{L} )^k \E[\dist(x_0,\domain\opt)]$, where $c > 0$ is a
numerical constant, which is the rate of convergence for (deterministic)
gradient methods with optimal stepsize choices~\cite{Nesterov04}.
More generally, we see a roughly linear speedup in the batch size $\mb$
to achieve a given accuracy until $\mb \ge \frac{\vargrowth^2}{\scvx \smooth}$:
to obtain $\E[\dist(x_k, \xdomain\opt)^2] \le \epsilon$ takes
\begin{equation*}
  k = O(1) \max\left\{ \frac{\lipgrad}{\lambda}, \frac{\vargrowth^2}{\lambda^2 \mb}
  \right\} \log \frac{1}{\epsilon}
\end{equation*}
iterations with appropriately chosen stepsize $\stepsize$. That is, we
expect to see a linear improvement in the number of iterations to achieve a
given accuracy $\epsilon$ until the condition number
$\frac{\lipgrad}{\lambda}$ dominates the variance of the gradient estimates.

\subsection{Upper bounds under an expected growth condition}
\label{sec:upper-expected-growth}

In Theorem~\ref{thm:easygammagrowth} above, we restrict the stepsizes
to have the form $\stepsize_k = \frac{1}{L + \eta_k}$. With more accurate
models and an alternative growth assumption on the functions $F$ and $f$, we can
remove this weakness, highlighting the robustness of more accurate models.
To that end, we revisit a few results of Asi and Duchi~\cite{AsiDu19siopt},
beginning with a slight generalization of their growth assumption
(which corresponds to the choices $\gamma \in \{0, 1\}$ below):
\begin{assumption}[$\gamma$-Growth]\label{assm:gamma-growth}
  There exist constants $\lambda_0,\lambda_1 > 0$ and $\gamma \in [0, 1]$,
  such that for all $\alpha \in \R_+, x \in \domain, x\opt \in \domain\opt$, we
  have
  \begin{equation*} 
    \begin{split}
      & \E\left[(F(x;\statrv) - F(x\opt;\statrv))\min\left\{\stepsize,\frac{F(x;\statrv) - F(x\opt,\statrv)}{\ltwo{F'(x;\statrv)}^2}\right\}\right]
      \\ 
      & \qquad\qquad ~ \geq \min \{\lambda_0\stepsize,\lambda_1{\rm dist}(x,\domain\opt)^{1-\gamma}\} \dist(x,\domain\opt)^{1 + \gamma}.
    \end{split}
  \end{equation*}
\end{assumption}

\noindent
As we will show in the coming section, while
Assumption~\ref{assm:gamma-growth} looks like a technical assumption, it
actually fairly closely governs the complexity of solving interpolation
problems, in that the $\lambda_1$ parameter describes lower bounds on the
convergence of any method. Essentially, the assumption states that the
functions $F$ must grow relative to the magnitude of their gradients at a
particular rate, so that it provides a type of stochastic growth condition.
We shall revisit this in the next section when we prove our lower bounds,
for now focusing on algorithms and their convergence under the
assumption. First, however, we may again rely on linear regression-type
objectives for an example satisfying Assumption~\ref{assm:gamma-growth}.

\begin{example}
  Consider a problem with data $\statval = (a, b) \in \R^n \times \R$,
  where $b = \<a, x\opt\>$ for all $(a, b)$, and set
  $F(x; (a, b)) = \frac{1}{1 + \gamma} |\<a, x - x\opt\>|^{1 + \gamma}$,
  so $\ltwo{F'(x; (a, b))}^2 = \ltwo{a}^2 |\<a, x - x\opt\>|^{2\gamma}$.
  If $a \sim \normal(0, I_n)$, then
  $|\<a, x - x\opt\>| \ge \half \ltwo{x - x\opt}$ with probability at least
  $\frac{3}{5}$, and similarly $\ltwo{a}^2 \le 2n$ with probability at least
  $\frac{3}{5}$, so that both occur with probability at least $\frac{1}{5}$.
  We then obtain
  \begin{align*}
    \E\left[F(x; \statrv) \min\left\{\stepsize,
      \frac{F(x; \statrv)}{\ltwo{F'(x; \statrv)}^2} \right\}\right]
    \ge \frac{1}{5} \frac{\ltwo{x - x\opt}^{1 + \gamma}}{
      2^{1 + \gamma} (1 + \gamma)}
    \min\left\{\stepsize, \frac{\ltwo{x - x\opt}^{1 - \gamma}}{
      2^{1 - \gamma} (1 + \gamma) \cdot 2n } \right\},
  \end{align*}
  so that Assumption~\ref{assm:gamma-growth} holds with
  $\lambda_0 \ge \frac{1}{5 (1 + \gamma) 2^{1 + \gamma}}$ and
  $\lambda_1 \ge \frac{1}{2^{2 - \gamma} (1 + \gamma) n}$.
\end{example}

To give stronger convergence results under Assumption~\ref{assm:gamma-growth},
we require one additional condition on our models, which
Asi and Duchi~\cite{AsiDu19siopt} introduce:
\begin{enumerate}[label=(C.\roman*),start = 3]
\item \label{cond:trunc} For all $\statval \in \statdomain$, the models
  $F_x(\cdot;\statval)$ satisfy
  \begin{equation*}
    F_x(y; \statval) \ge \inf_{z \in \domain} F(z; \statval).
  \end{equation*}
\end{enumerate}
In minibatch settings, where one considers a batch $\statrv^{1:\mb}$ of
samples in each model, the condition~\ref{cond:trunc} can be somewhat
challenging to verify, as it requires accuracy for the average $\inf_z
\wb{F}(z; \statval^{1:\mb})$, though (obviously) proximal methods
\eqref{eqn:full-prox} satisfy this condition, and in typical situations
(e.g.\ linear regression) where the batch size $\mb \le n$, the average of
truncated models~\eqref{eq:avg-model} will be similarly accurate.

% We extend a result by \cite{AsiDu19siopt} to form the following corollary, 
\begin{corollary}\label{corr:gamma-growth-convergence}
  Let Assumption~\ref{assm:gamma-growth} hold, and let $x_k$ be generated by
  the stochastic iteration~\eqref{eqn:model-iteration} for a model
  satisfying conditions \ref{cond:convex-model}--\ref{cond:trunc}. Take
  stepsizes $\alpha_k = \alpha_0 k^{-\beta}$ for some $\beta \in
  [0, 1]$. Define $K_0 \coloneqq \lfloor (\lambda_0\alpha_0/(\lambda_1
  \dist(x_1,\domain\opt)^{1 - \gamma}))^{1/\beta} \rfloor$.  Then
  \begin{align*}
    \E[\dist(x_{k+1},\domain\opt)^2] \le
    \exp{\left( - \lambda_1 \min\{k,K_0\} - \frac{\lambda_0}{\dist(x_1,\domain\opt)^{1 - \gamma}} \sum_{i = K_0 + 1}^{k}\alpha_i \right)}\dist(x_1,\domain\opt)^2.
\end{align*}
\end{corollary}
\begin{proof}
  Let $D_k = \dist(x_k,\domain\opt)$ and
  $\mc{F}_k = \sigma(\statrv_1, \ldots, \statrv_k)$ be the $\sigma$-field
  generated by the first $k$ samples $\statrv_i$. Then Lemma 4.1 of the
  paper~\cite{AsiDu19siopt} immediately yields
  \begin{align*}
    \E[D_{k+1}^2 \mid \mc{F}_{k-1}]
    &\le D_k^2 - \min\{\lambda_0 \stepsize_k D_k^{1+\gamma}, \lambda_1D_k^2\}.
  \end{align*}
  As $D_1 \ge D_k$ (again, by \cite{AsiDu19siopt}, Lemma 4.1), we in
  turn obtain
  \begin{align*}
    \E[D_{k + 1}^2 \mid \mc{F}_{k-1}]
    & \le \max\left\{1 - \lambda_1, 1 - \lambda_0 \stepsize_k / D_1^{1 - \gamma}
    \right\} D_k^2.
  \end{align*}
  The remainder of the argument is algebraic manipulations, as in the
  proof of Proposition 2 from \cite{AsiDu19siopt}.
\end{proof}

In the best case---when the stepsizes $\stepsize_k \uparrow \infty$ in
Corollary~\ref{corr:gamma-growth-convergence}---we achieve convergence
scaling as $\E[\dist(x_k, \mc{X}\opt)^2] \lesssim \exp(-\lambda_1 k)
\dist(x_1, \mc{X}\opt)^2$, and moreover (as we show in the next section)
this dependence on the growth constant $\lambda_1$ is unimprovable. With
this as motivation, one might hope that increased sampling (minibatching)
might increase the growth constant $\lambda_1$ in
Assumption~\ref{assm:gamma-growth}; here we provide a sketch of such a
result, which also makes it somewhat easier to check the conditions of
Assumption~\ref{assm:gamma-growth}, by giving three growth conditions.

\begin{enumerate}[label={(G.\roman*)}]
\item \label{cond:growth-sufficient1} There exists $\mu > 0$ and a
  probability $p > 0$ such that for all $x \in \domain$, we have
  \begin{equation*}
    \P(F(x; \statrv) - F(x\opt; \statrv) \geq \mu
    \dist(x, \domain\opt)^{1 + \gamma}) \geq p.
  \end{equation*}
\item \label{cond:growth-sufficient2} The (sub)gradient $f'$ is $(L,
  \gamma)$-Holder continuous, meaning $\ltwo{ f'(x) - f'(y)} \leq L \ltwo{x
    - y}^{\gamma}$, and $0 \in \partial f(x\opt)$.
\item \label{cond:growth-sufficient3} There exists $\noisetosig$ such
  that $\noisetosig \geq \frac{\var(F'(x;\statrv))}{\ltwo{f'(x)}^2}$ for all
  $x \in \domain$.
\end{enumerate}
\noindent
Our typical situation is to think of $\mu$ and $p$ numerical constants,
where the scaling $\noisetosig$ measures the noise inherent to the problem.
In any case, a short calculation shows how
Conditions~\ref{cond:growth-sufficient1}--\ref{cond:growth-sufficient3}
suffice to give Assumption~\ref{assm:gamma-growth}.

\begin{lemma}\label{lem:gamma-growth-ub}
  Let
  conditions~\ref{cond:growth-sufficient1}--\ref{cond:growth-sufficient3}
  hold. Then the average
  $\Fbar(x; \statval^{1:m}) = \frac{1}{m}\sum_{i = 1}^{n}F(x;\statval^i)$
  satisfies
  the $\gamma$-growth condition of Assumption~\ref{assm:gamma-growth}
  with
  \begin{equation*}
    \lambda_0 = \frac{\floor{mp}}{4m} \mu
    ~~ \mbox{and} ~~
    \lambda_1 = \frac{(\floor{mp} / m)^2 \mu^2}{16 \lipgrad^2 (1 + \frac{\noisetosig}{m}
      )}.
  \end{equation*}
  %%
  %% \begin{align}
  %%   \label{eq:sharp_growth_m_dep}
  %%   \E\left[\min\left\{\stepsize
  %%     \Fbar(x; \statrv^{1:m}), \frac{\Fbar(x; \statrv^{1:m})^2}{
  %%       \ltwobig{\Fbar'(x; \statrv^{1:m})}^2}\right\}
  %%     \right]
  %%   & \geq\frac{p}{2}
  %%   \min\left\{\stepsize \mu, \frac{p \mu^2 \dist(x, \domain\opt)^{1-\gamma}}{
  %%     2(1 + \frac{\noisetosig}{m})}\right\} \dist(x, \domain\opt)^{1+\gamma}.
  %% \end{align}
\end{lemma}
\begin{proof}
  For shorthand, we assume w.l.o.g.\ that $F(x\opt; \statrv) = 0$ with
  probability 1.
  The event that $F(x; \statrv^i) \ge \mu \dist(x, \xdomain\opt)^{1 + \gamma}$
  has probability at least $p$, and as the median
  of a $\binomdist(m,p)$ distribution lies in
  $\{\floor{mp}, \ceil{mp}\}$, we have
  \begin{equation}
    \P\left(\Fbar(x; \statrv^{1:m}) \ge \frac{\floor{mp}}{m}
    \mu \dist(x, x\opt)^{1 + \gamma}
    \right) \ge \half.
    \label{eqn:binomial-probs}
  \end{equation}
  Thus, the event
  \begin{equation*}
    A \defeq \left\{\ltwobig{\Fbar'(x;\statrv^{1:m})}^2 \leq
    4 \E\left[\ltwobig{\Fbar'(x; \statrv^{1:m})}^2\right],
    ~ \Fbar(x;\statrv^{1:m}) \geq \frac{\floor{mp}}{m} \mu
    \dist(x,x\opt)^{1+\gamma}
    \right\}
  \end{equation*}
  satisfies
  \begin{align*}
    \P(A)
    = 1 - \P(A^c)
    & \ge 1 - \P\left(\ltwobig{\Fbar'(x;\statrv^{1:m})}^2 \ge
    4 \E\left[\ltwobig{\Fbar'(x; \statrv^{1:m})}^2\right]\right)
    - \half
    \ge \frac{1}{4},
  \end{align*}
  where we use inequality~\eqref{eqn:binomial-probs}.
  We also have
  \begin{align*}
    \E\left[\ltwobig{\wb{F}'(x; \statrv^{1:m})}^2\right]
    = \ltwo{f'(x)}^2 \left(1 + \frac{\var(F'(x; \statrv))}{\ltwo{f'(x)}^2}
    \right)
    & \le \left(1 + \frac{\noisetosig}{m}\right) \ltwo{f'(x)}^2 \\
    & \le \left(1 + \frac{\noisetosig}{m} \right) \lipgrad^2
    \dist(x, \xdomain\opt)^{2\gamma},
  \end{align*}
  where we have used Conditions~\ref{cond:growth-sufficient3}
  and~\ref{cond:growth-sufficient2}.  Applying these
  observations gives
  \begin{align*}
    \lefteqn{\E\left[\min\left\{\stepsize
        \Fbar(x; \statrv^{1:m}), \frac{\Fbar(x; \statrv^{1:m})^2}{
          \ltwobig{\Fbar'(x; \statrv^{1:m})}^2}\right\}
        \right]} \\
    & \ge \frac{1}{4}
    \min\left\{ \stepsize \frac{\floor{mp}}{m}
    \mu \dist(x, \xdomain\opt)^{1 + \gamma},
    \frac{(\floor{mp} / m)^2 \mu^2 \dist(x, \xdomain\opt)^{2 + 2 \gamma}}{
      4 \E[\ltwos{\wb{F}'(x; \statrv^{1:m})}^2]} \right\} \\
    & \ge \frac{1}{4} \min \left\{
    \stepsize \frac{\floor{mp}}{m}
    \mu \dist(x, \xdomain\opt)^{1 + \gamma},
    \frac{(\floor{mp} / m)^2 \mu^2 \dist(x, \xdomain\opt)^2}{
      4 \lipgrad^2 (1 + \frac{\noisetosig}{m})}\right\},
  \end{align*}
  as desired.
  %% \begin{align*}
  %%   & \E\left[\min\left(\alpha \Fbar(x; \statrv_1^m), \frac{\Fbar(x; \statrv_1^m)^2}{ \ltwo{\Fbar'(x; \statrv_1^m)}^2}\right) \right] \\
  %%   %   \geq & \E\left[\min\left(\alpha \Fbar(x; \statrv_1^m), \frac{\Fbar(x; \statrv_1^m)^2}{ \ltwo{\Fbar'(x; \statrv_1^m)}^2}\right)\ind{A}\right] \\
  %%   \stackrel{(i)}\geq & \frac{p}{2} \min\left(\alpha c \dist(x, \domain\opt)^{1+\gamma}, \frac{pc^2 \dist(x, \domain\opt)^{2+2\gamma}}{2\ltwo{f'(x)}^2(1 + \frac{\var(F'(x;\statrv))}{m\ltwo{f'(x)}^2})}\right)\\
  %%   \stackrel{(ii)}\geq & \frac{p}{2} \min\left(\alpha c \dist(x, \domain\opt)^{1+\gamma}, \frac{pc^2 \dist(x, \domain\opt)^{2}}{2L^2(1 + \frac{\var(F'(x;\statrv))}{m\ltwo{f'(x)}^2})}\right)\\
  %%   \geq & \frac{p}{2} \min\left(\alpha c , \frac{pc^2 \dist(x, \domain\opt)^{1-\gamma}}{2L^2(1 + \frac{\noisetosig}{m})}\right)\dist(x, \domain\opt)^{1+\gamma},
  %% \end{align*}
  %% where (i) follows from $\E[X] \geq \E[X\ind{A}]$ and (ii) follows from Condition \ref{cond:growth-sufficient2}. \proofbox
\end{proof}

In brief, we see that mini-batches of size $\mb$ suggest improved
convergence related to the
noise-to-signal ratio $\noisetosig \defeq \sup_x \frac{\var(F'(x;
  \statrv))}{\ltwo{f'(x)}^2}$: once the sample size $m$ is large enough that
$\noisetosig / m \lesssim 1$, we expect relatively little improvement, though we
\emph{do} see a linear improvement in the growth constant $\lambda_1$ as $m$
grows whenever $m \ll \noisetosig$. To see this, let us for simplicity assume that
in Conditions~\ref{cond:growth-sufficient1}--\ref{cond:growth-sufficient3}
we have $p \gtrsim 1$ and $\lipgrad / \mu \lesssim 1$ (that is, the problem
is well-conditioned). Then applying
Corollary~\ref{corr:gamma-growth-convergence}, we see that
for large enough stepsizes $\stepsize$,
\begin{equation}
  \label{eqn:nice-conditioned-funny}
  k = O(1) \left(1 + \frac{\noisetosig}{m}\right) \log \frac{1}{\epsilon}
\end{equation}
iterations of any model-based method~\eqref{eqn:model-iteration} with
minibatches of size $m$---assuming that
Conditions~\ref{cond:convex-model}--\ref{cond:trunc} hold for the models
$\wb{F}_x$---are sufficient to guarantee $\E[\dist(x_k, \mc{X}\opt)^2]
\le \epsilon$.

\section{Optimality in Interpolation Problems}
\label{sec:lower-bounds}

We conclude the theoretical portion of this paper by developing several new
optimality results for interpolation problems, that is, those satisfying
Definition~\ref{def:interpolation}.  In brief, we shall show that the
depndence of Corollary~\ref{corr:gamma-growth-convergence} on the growth
constant $\lambda_1$ is sharp and unimprovable, and that in some cases, the
dependence on the signal-to-noise ratio $\noisetosig^{-1} \defeq \inf_x
\frac{\ltwos{f'(x)}^2}{\var(F'(x; \statrv))}$ is essentially sharp as well.
We do so via information-theoretic lower bounds on estimation of optimal
points, the first in a stylized $n = 1$ dimensional problem that gives the
correct dependence on the growth constants in
Assumption~\ref{assm:gamma-growth},
the second in standard regression problems but where
we choose the dimension $n \in \N$ more carefully.

We define our \emph{minimax risk} as follows.  Let
$\mc{P}$ be a family of problems, where a problem is a pair $(F, P)$
consisting of a probability distribution $P$ supported on $\statdomain$ and
function $F$ as defined in the introduction. We let $\domain\opt(F, P) =
\argmin_{x \in \domain} \E_P[F(x; \statrv)]$ be the collection of
minimizers, and define the minimax squared error
\begin{equation}
  \label{eqn:minimax-risk}
  \minimax_k(\mc{P}, \domain)
  \defeq \inf_{\what{x}^k} \sup_{(F, P) \in \mc{P}}
  \E_{P^k}\left[\dist(\what{x}^k, \domain\opt(F, P))^2\right],
\end{equation}
where the infimum is over all measurable $\what{x}^k :
\statdomain^k \to \R^n$, the supremum is over problems $(F, P) \in
\mc{P}$, and the inner expectation is over the samples $\statrv_1, \ldots,
\statrv_k \simiid P$.

\subsection{A lower bound for one-dimensional problems}

We first focus on problems for which we can isolate the contributions
of the growth constant $\lambda_1$ in Assumption~\ref{assm:gamma-growth},
letting the dimension $n = 1$ to show that our complexity bounds
hold independent of dimension; higher dimensions can only
yield increased complexity.
We consider a collection of well-conditioned problems, where we analogize
the typical condition number of $f$ by defining
\begin{equation*}
  \lambda_\gamma(f) \defeq \inf_{x \not \in \xdomain\opt} \frac{f(x) - f(x\opt)}{
    \frac{1}{1 + \gamma} \dist(x, \xdomain\opt)^{1 + \gamma}}
  ~~ \mbox{and} ~~
  \lipgrad_\gamma(f) \defeq \sup_{x \neq y} \frac{|f'(x) - f'(y)|}{|x - y|^\gamma},
\end{equation*}
calling $\kappa_\gamma(f) \defeq \frac{\lipgrad}{\lambda}$ the
condition number. We also note in passing that the constant $\lambda_1 \le
1$ in Assumption~\ref{assm:gamma-growth}, as by convexity we have
\begin{equation*}
  \frac{(F(x; \statval) - F(x\opt; \statval))^2}{F'(x; \statval)^2}
  \le \frac{\<F'(x; \statval), x - x\opt\>^2}{F'(x; \statval)^2}
  \le |x - x\opt|^2,
\end{equation*}
so taking $\stepsize \uparrow \infty$ in Assumption~\ref{assm:gamma-growth}
guarantees $\lambda_1 \in [0, 1]$.
Thus, for our first collection of problems, we let
\begin{equation}
  \label{eqn:nice-collection}
  \mc{P}_\gamma(\lambda_1)
\end{equation}
be those problems satisfying Assumption~\ref{assm:gamma-growth} with a given
$\gamma, \lambda_1 \in [0, 1]$,
any $\lambda_0 \ge \lambda_1$, our standing assumption of the
interpolation condition in Definition~\ref{def:interpolation},
and condition number $\kappa_\gamma(f) = 1$. The choice of the
condition number serves to highlight the difficulties from
stochasticity in the problem, eliminating the contributions of hardness
from the population (deterministic) objective $f$; an
identical lower bound will of course hold in the coming theorem for
more poorly conditioned problems with $\kappa_\gamma(f) \ge 1$,
as this is simply a larger collection.
\begin{theorem}
  \label{theorem:main-lower-bound}
  Let $\mc{P}_\gamma(\lambda_1)$ be the
  collection~\eqref{eqn:nice-collection}, assume that
  $\xdomain$ contains an $\ell_2$-ball of radius $\radius \ge 0$.
  Then
  \begin{equation*}
    \minimax_k(\mc{P}_\gamma(\lambda_1), \xdomain)
    \ge \frac{\radius^2}{2} \hinge{1 - (1 + \gamma)^2 \lambda_1}^k.
  \end{equation*}
\end{theorem}
\noindent
We make a few remarks before proceeding to the proof. First, the convergence
guarantees in Section~\ref{sec:upper-expected-growth} show that appropriate
model-based methods converge to $\epsilon$ accuracy in
$O(\frac{1}{\lambda_1} \log \frac{1}{\epsilon})$ iterations, which
by the theorem is optimal. Thus, in a strong sense, the \emph{a priori}
esoteric-seeming growth condition in Assumption~\ref{assm:gamma-growth}
is indeed fundamental.

\begin{proof}
  Let $\mc{P} = \mc{P}_\gamma(\lambda_1)$ for short, and
  assume w.l.o.g.\ that $\lambda_1 \le 1/(1 + \gamma)^2$, as the result is
  trivial otherwise. We base our argument on
  Le Cam's two point method (see, e.g., \cite{Wainwright19}, Eq.~(15.14)). We
  consider two probability distributions $P_1, P_{-1}$, and let
  $\mc{X}_v\opt$ be (for now) arbitrary sets indexed by $v \in \{\pm 1\}$.
  Then recall the variation distance $\tvnorm{P - Q} = \sup_A |P(A) -
  Q(A)|$ between distributions $P$ and $Q$, we have
  Le Cam's two-point method:
  \begin{lemma}[Le Cam]
    \label{lemma:le-cam}
    Let $\what{x}^k$ be an arbitrary function
    of $\statrv_1, \ldots, \statrv_k$.
    Then
    \begin{equation*}
      \max_{v \in \{-1, 1\}}
      \E_{P_v^k}\left[\dist(\what{x}^k, \domain_v\opt)^p\right]
      \ge \frac{1}{8}
      \dist(\domain_{-1}\opt, \domain_1\opt)^2
      \left(1 - \tvnorms{P_{-1}^k - P_1^k}\right).
    \end{equation*}
  \end{lemma}

  To use Lemma~\ref{lemma:le-cam} to lower bound the minimax risk
  it suffices to choose a pair of problems $(F, P_v) \in
  \mc{P}$ whose optimal sets are well-separated and
  apply the lemma. To that end, let $\delta \in (0, 1)$
  to be chosen later, and consider the choices
  \begin{equation}
    \label{eqn:near-point-distributions}
    P_{-1} :
    \begin{cases} \statrv = 0 & \mbox{w.p.}~ 1 - \delta \\
      \statrv = -1 & \mbox{w.p.}~ \delta
    \end{cases}
    ~~~~
    P_1: \begin{cases} \statrv = 0 & \mbox{w.p.}~ 1 - \delta \\
      \statrv = 1 & \mbox{w.p.}~ \delta.
    \end{cases}
  \end{equation}
  Our functions $F$ are trivial to construct: given the radius $\radius$,
  we define
  \begin{equation}
    \label{eqn:sharp-growth-lb-F}
    F(x; 1) = \frac{1}{1 + \gamma} |x - \radius|^{1 + \gamma}, ~~
    F(x; -1) = \frac{1}{1 + \gamma} |x + \radius|^{1 + \gamma}, ~~
    F(x; 0) = 0.
  \end{equation}
  The intuition here is that given a sample $\statrv \in \{-1, 0, 1\}$, we
  either completely identify the distribution or receive no information.

  It remains to show that the pairs $(F, P_v) \in \mc{P}$ and to
  bound the variation distance $\tvnorms{P_1^k - P_{-1}^k}$. For
  the latter, we have
  \begin{lemma}
    \label{lemma:exponential-tv}
    Let $P_{-1}, P_1$ be as in Eq.~\eqref{eqn:near-point-distributions}. Then
    $\tvnorms{P_{-1}^k - P_1^k} = 1 - (1 - \delta)^k$.
  \end{lemma}
  \begin{proof}
    For any distributions $P, Q$, with densities $p, q$ w.r.t.\ a base
    measure $\mu$, we have
    $\tvnorm{P - Q} = P(p > q) - Q(p > q)$.
    For $P_{-1}, P_1$ as above, we thus have
    \begin{align*}
      \tvnorms{P_{-1}^k - P_1^k}
      & = P_1^k(\mbox{there~exists}~i \in [k]
      ~ \mbox{s.t.}~ S_i = 1) \\
      & = 1 - P_1(S_1 = 0, \ldots, S_k = 0)
      = 1 - (1 - \delta)^k.
    \end{align*}
  \end{proof}

  Now, consider the functions
  \begin{equation*}
    f_v(x) \defeq \E_{P_v}[F(x; \statrv)]
    = \frac{\delta }{1 + \gamma} |x - v \radius|^{1 + \gamma}.
  \end{equation*}
  We have $\kappa_\gamma(f) = 1$, so that the problem is well-conditioned,
  and the optimal sets $\xdomain_v\opt \defeq \argmin_{x \in \xdomain} f_v(x)$
  are the singletons $\xdomain_v\opt = \{x_v\opt = v
  \radius\}$. Additionally, we have
  \begin{align*}
    \lefteqn{\E_v\left[(F(x; \statrv) - F(x_v\opt; \statrv))
        \min\left\{ \stepsize,
        \frac{F(x; \statrv) - F(x_v\opt; \statrv)}{
          \ltwo{F'(x; \statrv)}^2}\right\}\right]} \\
    & = \frac{\delta}{1 + \gamma} |x - v \radius|^{1 + \gamma}
    \min \left\{\stepsize, \frac{|x-v \radius|^{1+\gamma}}{(1 + \gamma)
      |x - vr|^{2\gamma}} \right\} \\
    & = \min\left\{\frac{\delta \stepsize}{1 + \gamma},
    \frac{\delta}{(1 + \gamma)^2}
    \dist(x, \domain_v\opt)^{1-\gamma} \right\}
    \dist(x, \domain_v\opt)^{1+\gamma},
  \end{align*}
  so by choosing $\delta = (1 + \gamma)^2 \lambda_1 \le 1$, our
  problems problems $(F, P_v)$ belong to $\mc{P}_\gamma(\lambda_1)$.
  Le Cam's Lemma~\ref{lemma:le-cam} and the variation
  distance bound in Lemma~\ref{lemma:exponential-tv} imply that
  \begin{equation*}
    \max_{v \in \{\pm 1\}}
    \E_{P_v^k}\left[|\what{x}^k - x_v\opt|^2\right]
    \ge \frac{1}{8} |x_1\opt - x_{-1}\opt|^2 (1 - \delta)^k
    = \frac{\radius^2}{2}(1 - \delta)^k.
  \end{equation*}
  Substituting $\delta = (1 + \gamma)^2\lambda_1$ gives the result.
\end{proof}

\subsection{A lower bound for well-conditioned regression problems}

\newcommand{\prior}{\pi}

The proof of Theorem~\ref{theorem:main-lower-bound} relies on constructing
certain power functions and a very careful choice of growth and probability.
An alternative approach is to mimic those ideas in proving complexity
results for deterministic problems~\cite{NemirovskiYu83, Nesterov04,
  CarmonDuHiSi19}, where one takes the dimension larger. By allowing
high-dimensional problems, we can show that the noise-to-signal ratio
$\noisetosig \coloneqq \sup_x \frac{\var(F'(x; \statrv))}{\norm{\nabla f(x)}^2}$
and growth constant $\lambda_1$ from Assumption~\ref{assm:gamma-growth} remain
fundamental, even in noiseless linear regression.

To make the proof cleaner we make a slight modification to the class of
problems we consider: instead of assuming a bounded domain $\xdomain$, we
instead assume $\xdomain = \R^n$, but now we consider a randomized (instead
of minimax/worst case) adversary that chooses a problem $(F, P) \in
\mc{P}$ according to a measure $\prior$ on the space of problems; in
particular, we assume that $\E_\mu[\ltwo{x_0 - x\opt}^2] \le \radius^2$,
that is, the expected distance of $x_0$ to $x\opt$ is at most $\radius$.
Letting $\xdomain\opt(F, P) = \argmin_x \E_P[F(x; \statrv)]$ be the optimal
set for a given problem $(F, P)$, we then define the minimum average risk
\begin{equation*}
  \minimax_k(\mc{P}, \prior) \defeq
  \inf_{\what{x}^k}
  \int \E_{P^k}[\dist(\what{x}^k, \xdomain\opt(F, P))^2] d\prior(F, P).
\end{equation*}
We note that the minimum average risk defined here naturally lower bounds the minimax risk ~\eqref{eqn:minimax-risk}, redefined analogously for our problem.

We specialize this randomized risk
for each $n \in \N$, letting
$\mc{P}_n$ be a collection of noiseless linear regression problems on
$\R^n$, where we identify the prior measure $\prior$
with $x\opt \sim \normal(0, \frac{\radius^2}{n} I_{n \times n})$.
Then certainly $\E[\ltwo{x\opt}^2] = \radius^2$.
We consider samples $\statval$ consisting of a pair $A \in \R^{\mb \times n}$
and $b = Ax\opt$, considering
the quadratic loss
\begin{equation}
  \label{eqn:squared-error}
  F(x; \statval) = F(x; (A, b)) = \half \ltwo{Ax - b}^2,
\end{equation}
and we call the resulting objective $f(x) = \E[F(x; \statrv)]$
\emph{perfectly conditioned} if $f(x) = c \ltwo{x - x\opt}^2$ for
a constant $c \in \R_+$.
We have the following theorem.

\begin{theorem}
  Let $\lambda_1 \in [0, \frac{1}{4}]$ and $\gamma = 1$. Then there exists a
  collection $\mc{P}$ of perfectly conditioned interpolating
  problems with squared error~\eqref{eqn:squared-error},
  satisfying Assumption~\ref{assm:gamma-growth}
  and $\E_\prior[\ltwo{x\opt}^2] = \radius^2$,
  such that
  \begin{equation*}
    \minimax_k(\mc{P}, \prior) \ge \radius^2 (1 - 4 \lambda_1)^k.
  \end{equation*}
  Alternatively, let $\noisetosig \in [1, \infty]$.  There exists a
  collection $\mc{P}$ of perfectly conditioned interpolating problems with
  squared error~\eqref{eqn:squared-error}, with
  noise-to-signal ratio satisfying $\sup_x \frac{\var(\grad F(x; \statrv))}{
    \ltwo{\grad f(x)}^2} \le \noisetosig$, such that
  \begin{equation*}
    \minimax_k(\mc{P}, \prior) \ge \radius^2 \left(1 - \frac{1}{\noisetosig}
    \right)^k.
  \end{equation*}
\end{theorem}
\noindent
Thus, one cannot hope
to achieve (much) better convergence even for quadratics than that
we have outlined: the dependence on either the growth
$\lambda_1$ or the signal-to-noise $\noisetosig^{-1}$ is unavoidable, and
one must collect at least $k \gtrsim \frac{1}{\lambda_1}
\log \frac{1}{\epsilon}$ or
$k \gtrsim \noisetosig \log \frac{1}{\epsilon}$ samples $\statrv$ to
achieve accuracy $\epsilon$, again
highlighting that these quantities---as we (inspired by
Asi and Duchi~\cite{AsiDu19siopt}) identify
in Corollary~\ref{corr:gamma-growth-convergence}
and the iteration bound~\eqref{eqn:nice-conditioned-funny}---are
fundamental for interpolation problems.

\begin{proof}
  Let $U = [u_1 ~ \cdots ~ u_n] \in \R^{n \times n}$ be an arbitrary orthogonal
  matrix, so $U^T U = UU^T = I_n$.  
  Let $\mc{P}_n$ be the collection of linear regression problems with data
  matrices $A \in \R^{\mb \times n}$ chosen by taking
  $\mb \le n$ columns $(u_{i(1)}, \ldots, u_{i(\mb)})$ of $U$ uniformly at random
  and setting $A = \sqrt{n/m} [u_{i(1)} ~ \cdots ~ u_{i(\mb)}]^T$,
  so that $\E[A^T A] = I_n$ and $(A^T A)^2 = (n/m) A^T A$,
  and
  let $b = A x\opt$, where
  $x\opt \sim \prior = \normal(0, \frac{\radius^2}{n} I_n)$ follows a Gaussian
  prior. Each observation $\statrv_i$ corresponds to releasing (perfectly) a
  random linear projection of $x\opt$, so that given the $k$
  observations,
  if we let $C_k = [A_1 ~ \cdots ~ A_k] \in \R^{n \times \mb k}$ denote the
  concatenated data matrix after $k$ observations,
  the posterior on $x\opt$ is
  \begin{equation*}
    x\opt \mid (\statrv_1, \ldots, \statrv_k)
    \sim \normal\left(\E[x\opt \mid \statrv_1, \ldots, \statrv_k],
    \frac{\radius^2}{n} (I_n - C_k (C_k^T C_k)^{-1} C_k^T)
    \right),
  \end{equation*}
  that is, the covariance projects out $C_k$.
  By a standard Bayesian argument,
  \begin{equation}
    \inf_{\what{x}^k}
    \E\left[\ltwobig{\what{x}^k - x\opt}^2 \right]
    = \E\left[\ltwobig{\E[x\opt \mid \statrv_1^k] - x\opt}^2\right]
    = \radius^2 \E\left[\frac{n - \rank(C_k)}{n}\right],
    \label{eqn:first-rank}
  \end{equation}
  as $I_n - C_k (C_k^T C_k)^{-1} C_k^T$ is a rank $n - \rank(C_k)$
  projection matrix.  
  Let $r_k = \rank(C_k)$ for shorthand. Then we may compute
  $\E[r_k]$ exactly by noting that
  \begin{equation*}
    \E[r_k \mid r_{k-1}]
    = r_{k-1} + \mb \frac{n - r_{k-1}}{n}
    = \left(1 - \frac{\mb}{n}\right) r_{k-1} + \mb,
  \end{equation*}
  so that with $r_1 = \mb$ we obtain
  \begin{equation*}
    \E[r_k] = \mb \sum_{i = 1}^k \left(1 - \frac{\mb}{n}\right)^{k - i}
    = \mb \frac{1 - (1 - \mb/n)^k}{1 - (1 - \mb/n)} =
    n - n\left(1 - \frac{\mb}{n}\right)^k,
  \end{equation*}
  and substituting this into expression~\eqref{eqn:first-rank} gives
  \begin{equation}
    \label{eqn:second-rank}
    \inf_{\what{x}^k} \E\left[\ltwobig{\what{x}^k - x\opt}^2 \right]
    = \radius^2 \left(1 - \frac{\mb}{n} \right)^k.
  \end{equation}

  We now use expression~\eqref{eqn:second-rank} to prove the two results in
  the theorem.
  For the first, we note that for $\statval = (A, b)$, we have
  $\grad F(x; \statval) = A^T (A x - b) = A^T A(x - x\opt)$, and as
  $(A^TA)^2 = \frac{n}{m} A^TA$ by construction and $\E[A^TA] = I_n$,
  \begin{align*}
    \lefteqn{\E\left[(F(x; \statrv) - F(x\opt; \statrv))
        \min\left\{\stepsize, \frac{F(x; \statrv) - F(x\opt; \statrv)}{
          \ltwo{\grad F(x; \statrv)}^2}\right\}\right]} \\
    & = \E\left[\min\left\{\frac{\stepsize}{2}
      \ltwo{A (x - x\opt)}^2,
      \frac{\ltwo{A(x - x\opt)}^4}{4 \ltwo{A^T A(x - x\opt)}^2}
        \right\}\right]
    = \min\left\{\frac{\stepsize}{2}, \frac{\mb}{4n}\right\}
    \ltwo{x - x\opt}^2.
  \end{align*}
  In particular,
  we can choose $m, n$ so that $\frac{\mb}{4n} \ge \lambda_1$
  the problem satisfies Assumption~\ref{assm:gamma-growth} with
  $\gamma = 1$ and $\lambda_0 = \half$. This gives the first result
  by substituting into expression~\eqref{eqn:second-rank}
  and taking $\mb, n$ so that $\frac{\mb}{n}$ is arbitrarily close
  to $4 \lambda_1$.

  For the second result,
  we recognize the noise-to-signal ratio
  \begin{equation*}
    \frac{\var(\grad F(x; \statrv))}{\ltwo{\grad f(x)}^2}
    \le \frac{\frac{n}{\mb} \ltwo{x - x\opt}^2}{\ltwo{x - x\opt}^2}
    = \frac{n}{\mb}.
  \end{equation*}
  Making appropriate substitutions by taking $\frac{n}{\mb} \le \rho$
  gives the second lower bound.
\end{proof}

\section{Experiments}
\label{sec:experiments}

\newcommand{\method}{\mathsf{a}}

Our goal now is to study and demonstrate the speedup and robustness of
\aprox methods with minibatches, comparing the relative performance of the
proposed methods on several benchmark stochastic optimization problems.
We consider the following five methods in our experiments, where
we use both single sample ($\mb = 1$) and minibatch ($\mb > 1$) versions:
\begin{enumerate}[leftmargin=2em]
\item SGM: stochastic gradient methods, i.e.,
  the linear model~\eqref{eqn:dumb-linear-model}.
\item Proximal: full proximal model~\eqref{eqn:full-prox} with
  averaged function \eqref{eq:avg}.
\item \PIA: truncated model~\eqref{eqn:trunc-model}
  with naive iterate averaging~\eqref{eq:iterate-avg}.
  %\eqref{eq:iterate-avg} with a truncated model \eqref{eqn:trunc-model} (avg-iter).
\item \PMA: iterates via truncating
  the averaged linear model, update~\eqref{eq:model-avg-trunc}.
  %Averaging the model \eqref{eq:avg-model} with a truncated model \eqref{eqn:trunc-model} (avg-model).
\item \PAM: iterates defined an average of truncated
  models, using updates~\eqref{eq:avg-model}.
  %Modeling the average \eqref{eq:model-avg} with a truncated model~\eqref{eqn:trunc-model} (model-avg).
\end{enumerate}

For our experiments, we use stepsizes $\stepsize_k = \stepsize_0 k^{-1/2}$,
varying $\stepsize_0$, and for each algorithm $\method$ report the number
$T_{\method, \mb}(\stepsize_0)$ of total samples consumed---as a proxy for
time---to reach $\varepsilon$ accuracy using minibatches of size $\mb$; that
is, $T_{\method,\mb}(\stepsize_0) = k\mb$ where $k$ is the first iteration
to satisfy $f(x_k) - f(x\opt) \le \varepsilon$. We also let $T\opt_{\method,
  \mb} = \min_{\stepsize_0} T_{\method, \mb}(\alpha_0)$ denote the smallest
time to convergence for a method $\method$ using batch size $\mb$.
Each of our experiments involves data $(A, b) \in \R^{N \times n} \times \R^N$,
where $f_{A,b}(x) = \frac{1}{N} \sum_{i = 1}^N F(x; a_i, b_i)$ for
a given loss $F$, and we vary the condition number of $A$,
taking $N = 10^3$ and $n = 40$.
With
these values identified, we present three types of results, focusing on
results that allow a more careful accouning for the robustness of the
various methods:

\begin{enumerate}[leftmargin=2em]
\item {\bf Performance profiles}~\cite{DolanMo02}\textbf{:} For each method
  $\method$, we evaluate for each $r \ge 1$ the fraction of the total
  executed experiments for which the $T_{\method, \mb}(\alpha_0) \leq r
  T_{\method\opt, \mb}(\alpha_0)$, where $\method\opt$ is the best
  performing method in each experiment, giving $r$ on the horizontal axis
  and the proportion on the vertical. Here, to evaluate robustness, we
  define a single experiment as one execution of each of the 5 methods for a
  particular step size $\stepsize_0$, minibatch size $\mb$, and condition
  number combination. We discard the experiments where more than 3 of the
  methods fail to complete before the max number of iterations.
\item {\bf Best speedups for minibatching:} For each method $\method$, we
  plot ${T\opt_{\method, 1}}/{T\opt_{\method,\mb}}$ against the
  minibatch size $\mb$ to show the speedup minibatching provides using the
  best step sizes. This shows the best possible speedup obtained by
  minibatching through tuning the initial step size $\stepsize_0$.
\item {\bf Time to solution \wrt~step-size:} For each method $\method$ and
  minibatch size $\mb$, we plot $T_{\method, \mb}(\stepsize_0)$ against the
  initial step size $\stepsize_0$.
\end{enumerate}

% {\bf Data Generation Methodology:}
We use minibatch sizes $\mb \in \{1,4,8,16,32,64\}$ and initial steps
$\stepsize_0 \in \{10^{i/2}, i \in \{-4, -3, \ldots, 5\}\}$. For all
experiments we run 30 trials with different seeds and plot the $95\%$
confidence sets. We describe the objective function and noise
mechanism for each problem in the respective subsections.

\subsection{Linear Regression}
\begin{figure}[ht]
  \begin{center}
    \begin{tabular}{cccc}
      \begin{overpic}[width=.45\columnwidth]{%,grid]{%
      		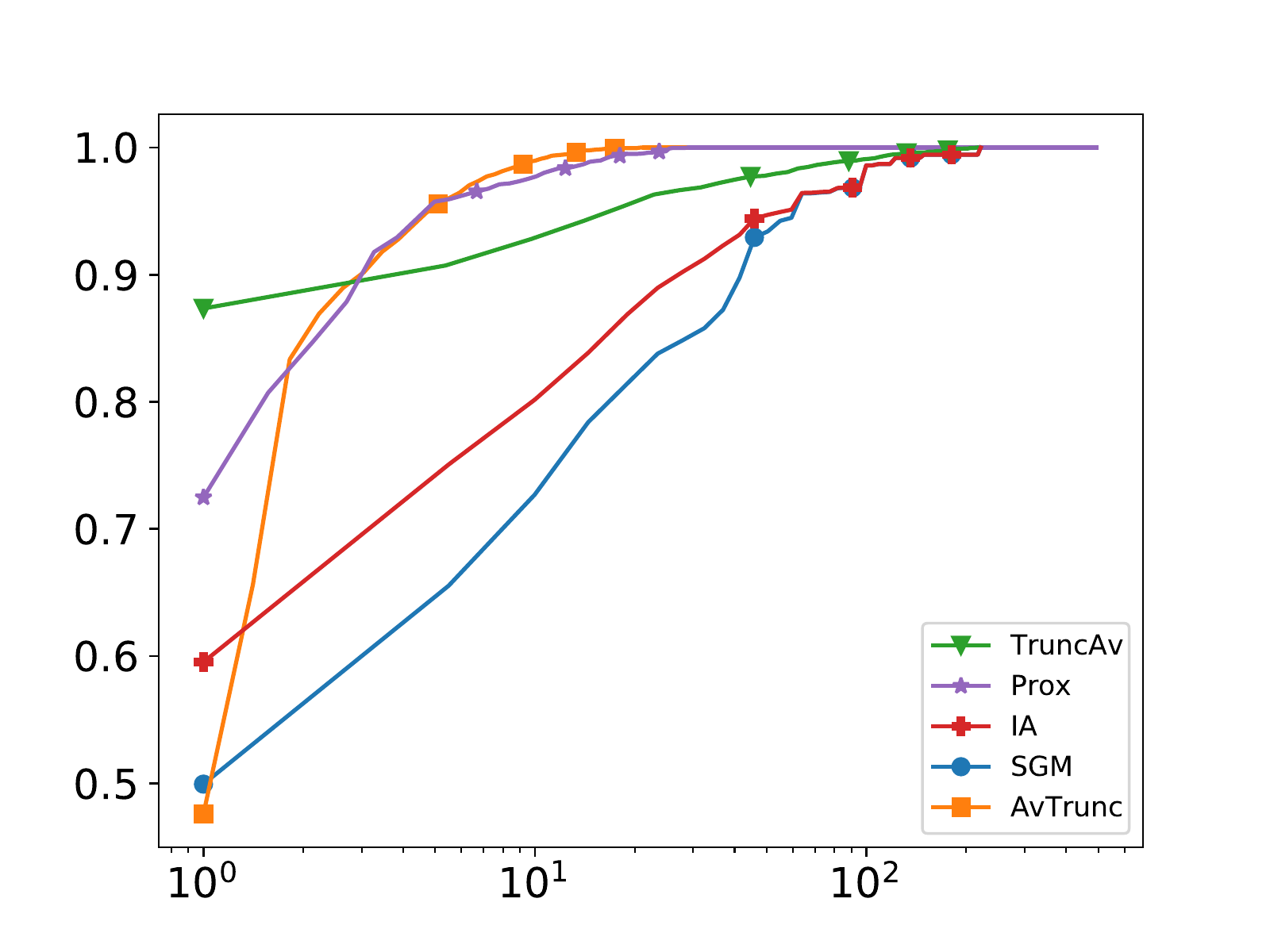}
      		\put(-2,8){
          \rotatebox{90}{{\small Fraction of Experiments}}}
        \put(28,-1){{\small Performance ratio $r$}}
        \put(32,68){{Non-accelerated}}
      \end{overpic} &
      \begin{overpic}[width=.45\columnwidth]{%,grid]{%
      		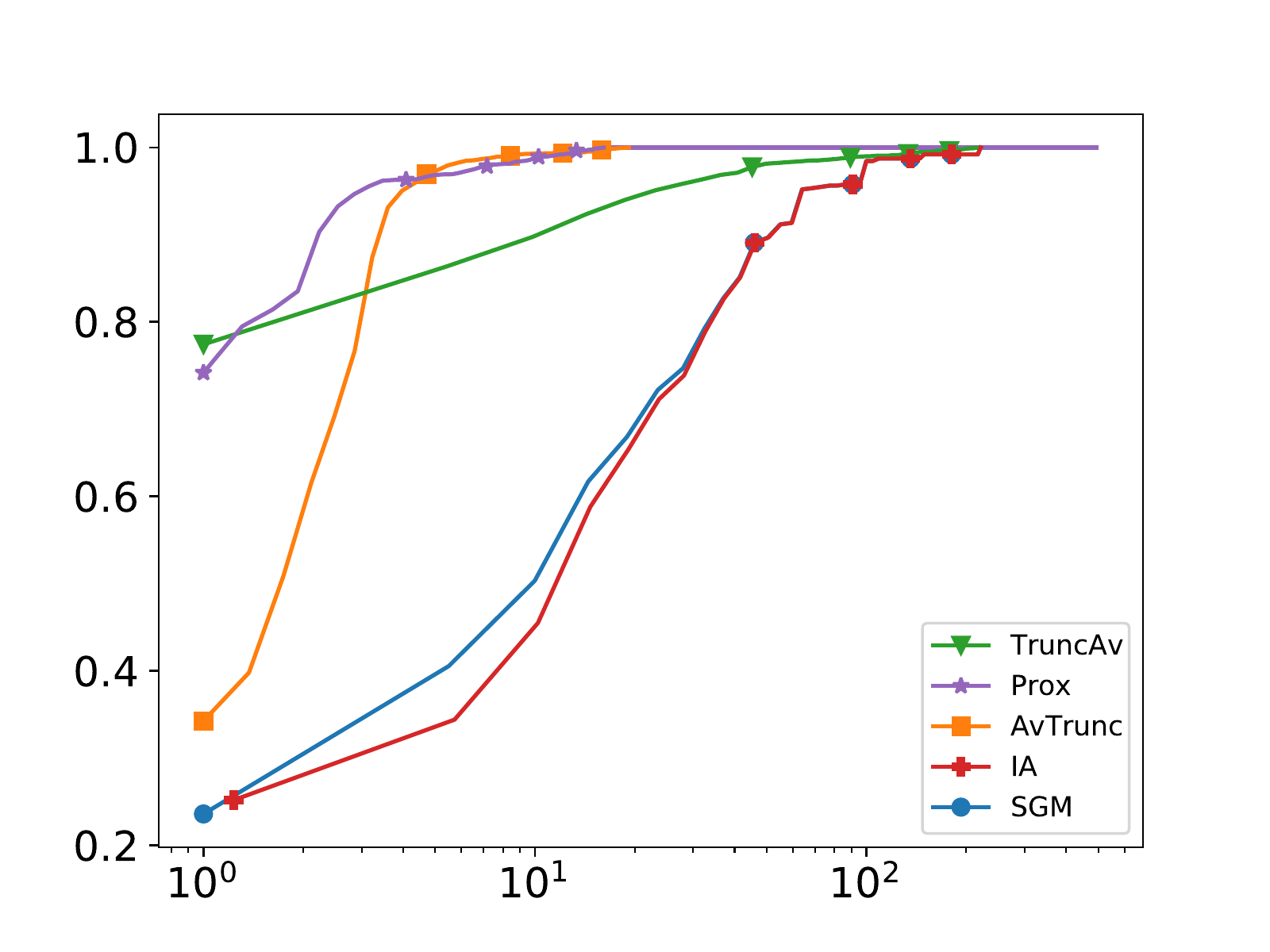}
      		\put(-2,8){
          \rotatebox{90}{{\small Fraction of Experiments}}}
          \put(28,-1){{\small Performance ratio $r$}}
          \put(37,68){{Accelerated}}
      \end{overpic} 
    \end{tabular}
  \caption{
    \label{fig:lin-perf-plot}
   Performance profiles for linear regression. 
    %for (a) noiseless absolute regression, (b) noiseless logistic regression, (c) noisy absolute regression and (d) noisy logistic regression.
  }
  \end{center}
  \vspace{-0cm}
\end{figure}

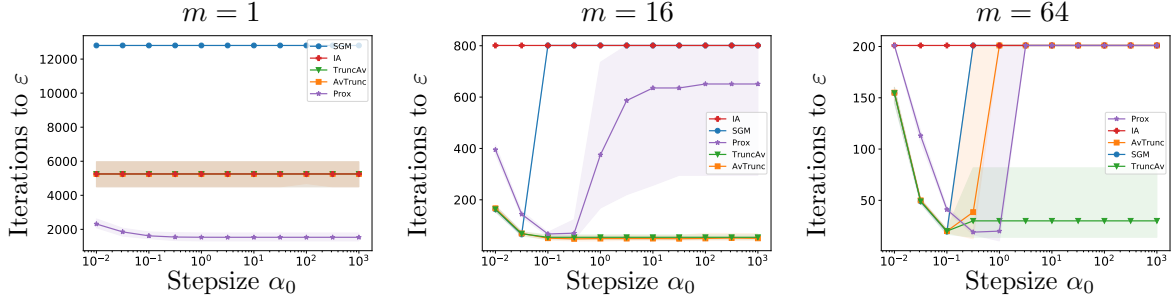
\begin{figure}[ht]
  \begin{center}
    \begin{tabular}{ccc}
      \begin{overpic}[width=.3\columnwidth]{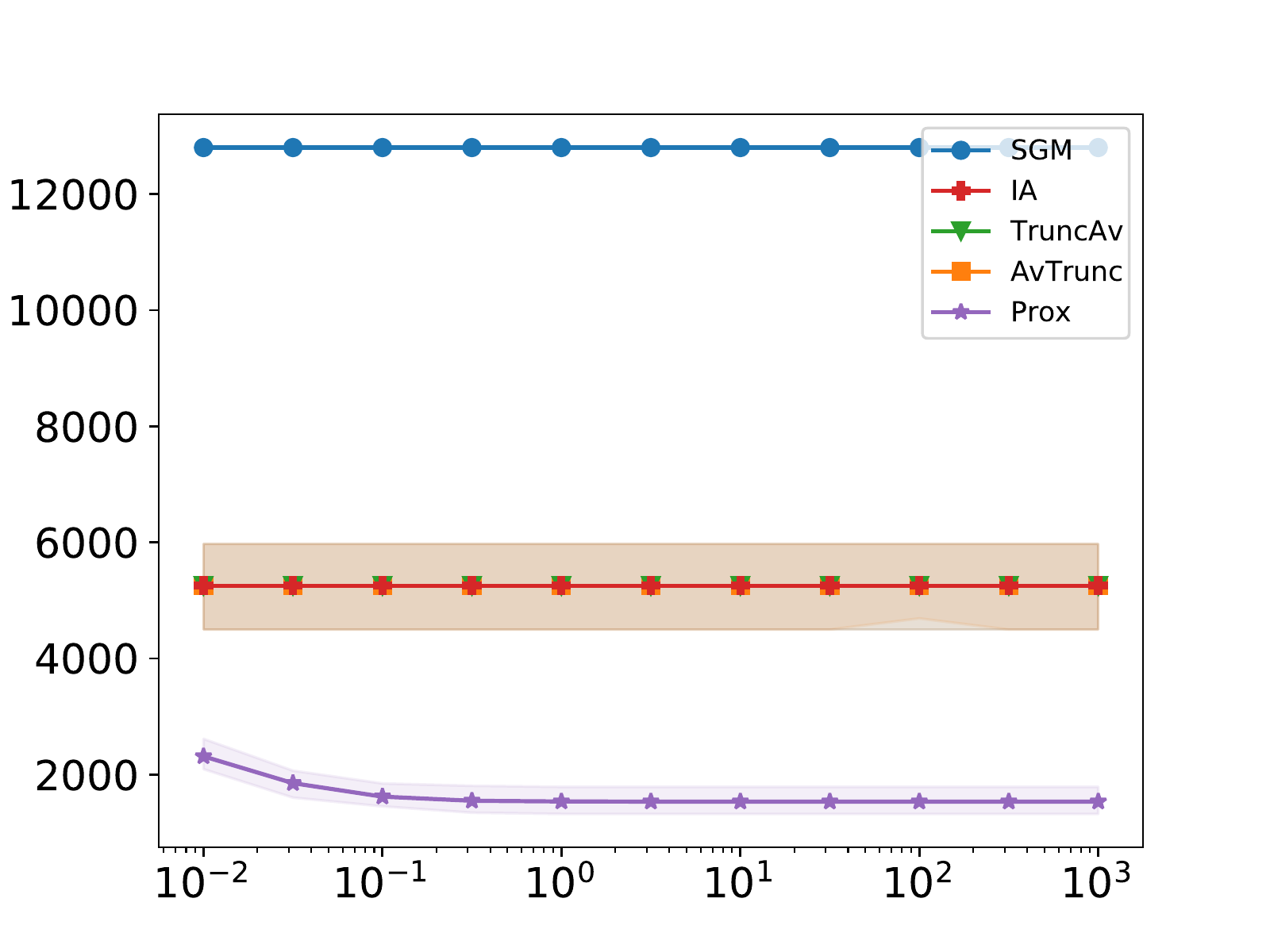}
      		\put(-10,12){
          \rotatebox{90}{{\small Iterations to $\varepsilon$}}}
      		\put(35,-2){{\small Stepsize $\stepsize_0$}}
      		\put(39,70){{$m=1$}}
      \end{overpic} &
      \begin{overpic}[width=.3\columnwidth]{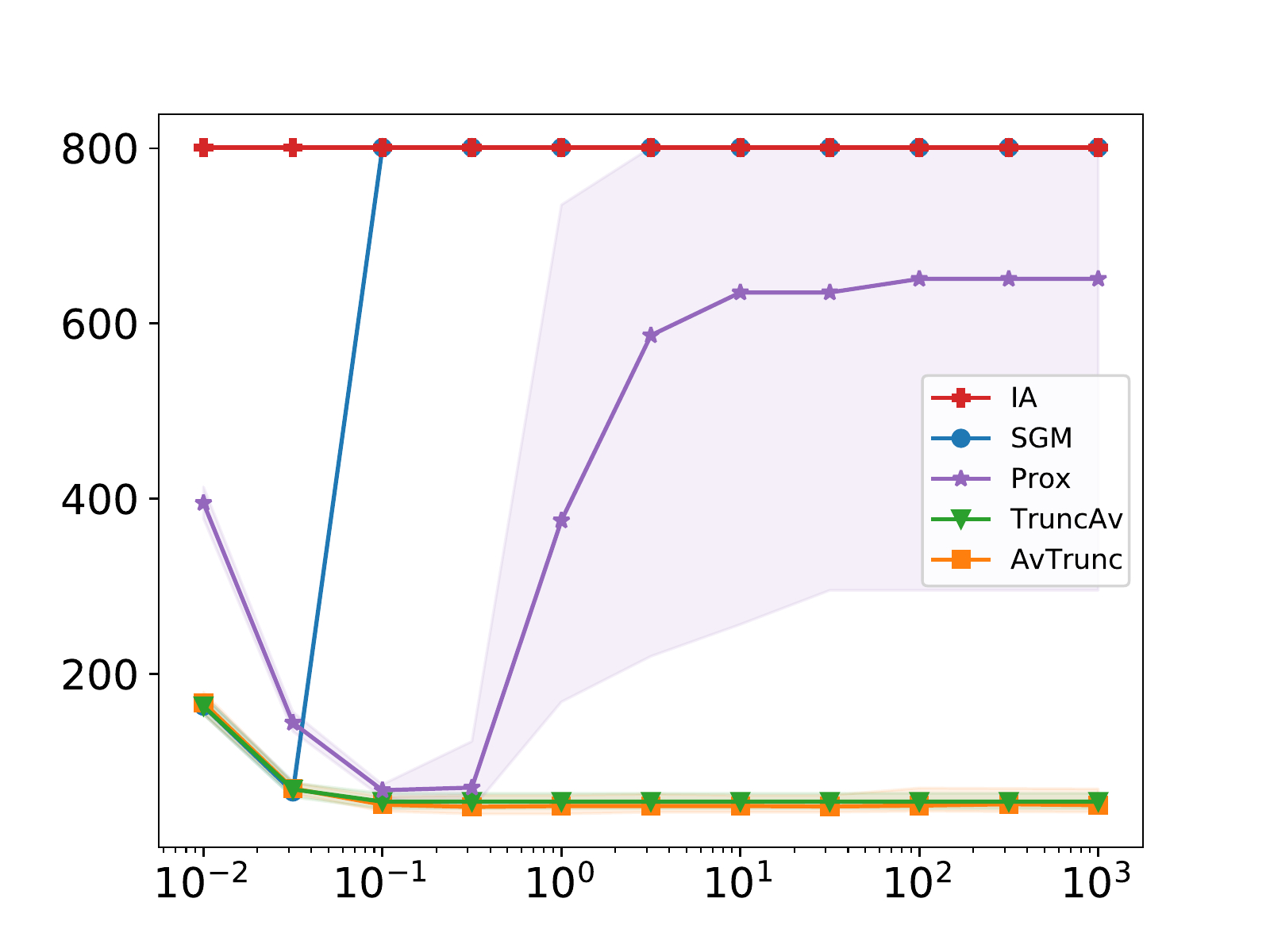}
      		\put(-10,12){
          \rotatebox{90}{{\small Iterations to $\varepsilon$}}}
      		\put(35,-2){{\small Stepsize $\stepsize_0$}}
      		\put(39,70){{$m=16$}}
      \end{overpic} &
      \begin{overpic}[width=.3\columnwidth]{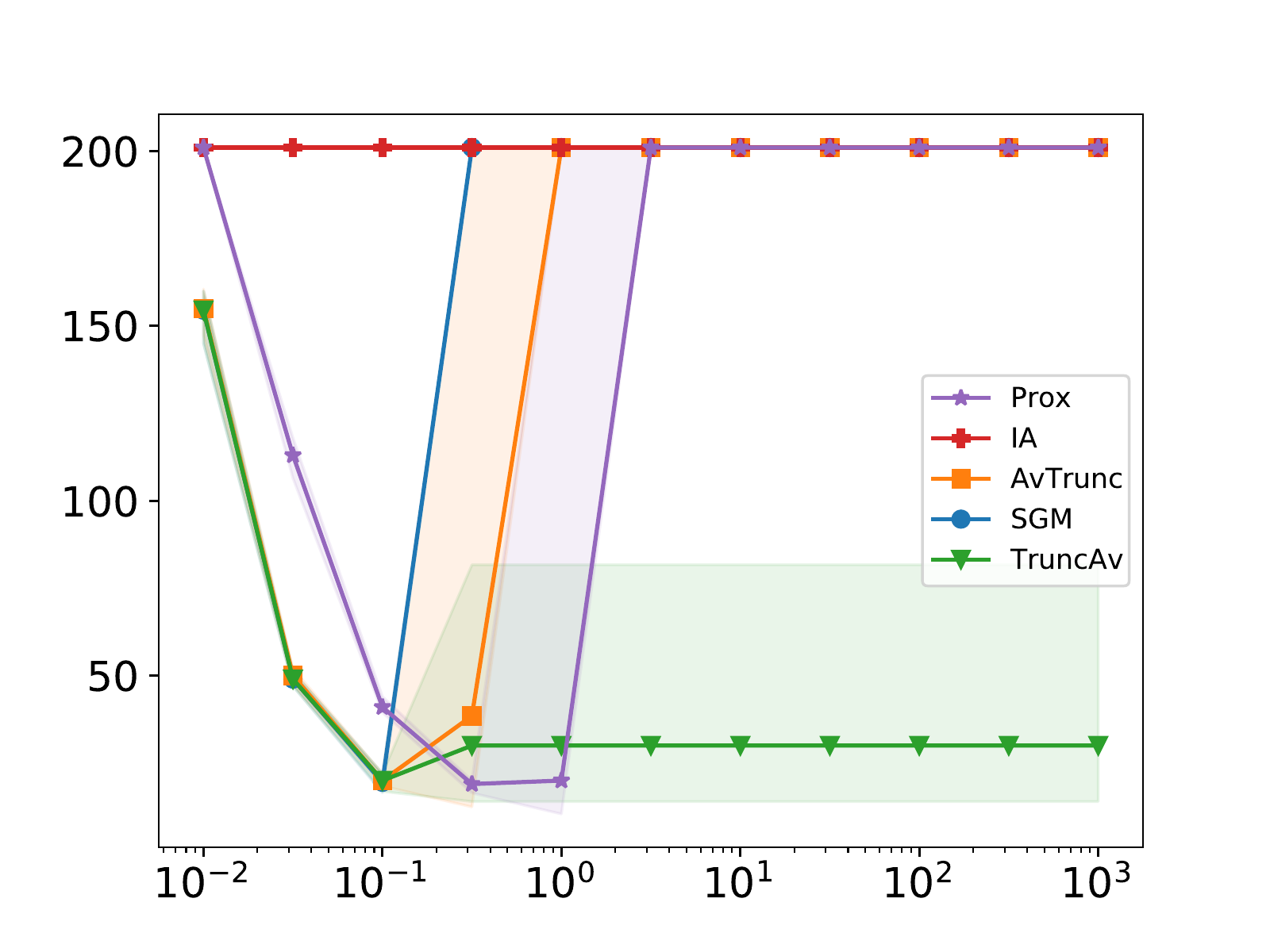}
      		\put(-10,12){
          \rotatebox{90}{{\small Iterations to $\varepsilon$}}}
      		\put(35,-2){{\small Stepsize $\stepsize_0$}}
      		\put(38,70){{$m=64$}}
      \end{overpic}
    \end{tabular}
  \caption{
    \label{fig:noiseless-linear-iteration}
    Time to convergence of the accelerated methods
    vs.\ stepsizes 
    for noisy linear regression 
    %. (a) $\mb = 1$, (b) $\mb = 8$, (c) $\mb = 32$
    }
  \end{center}
  \vspace{-.3cm}
\end{figure}

We have $f(x) = \frac{1}{2 N}\ltwo{Ax - b}^2$. For each experiment we
generate rows of $A$ and $x\opt$ i.i.d.\ $\normal(0,I_n)$ and, setting $b =
Ax\opt + \sigma v$ with $v \sim \normal(0, I_N)$. In the noisy setting for
our experiments, we set $\sigma = 0.5$. \Cref{fig:lin-perf-plot} outlines
the performance profiles for the linear regression experiments. The fully
proximal, \PAM, and \PMA methods are noticeably better than \PIA~and SGM.
\cref{fig:noiseless-linear-iteration} also reflects this behavior, where the
accelerated fully proximal, \PAM, and \PMA methods are more robust to
initial step size choice.

\subsection{Absolute loss regression}
\begin{figure}[ht]
  \begin{center}
    \begin{tabular}{cccc}
      \begin{overpic}[width=.45\columnwidth]{%,grid]{%
      		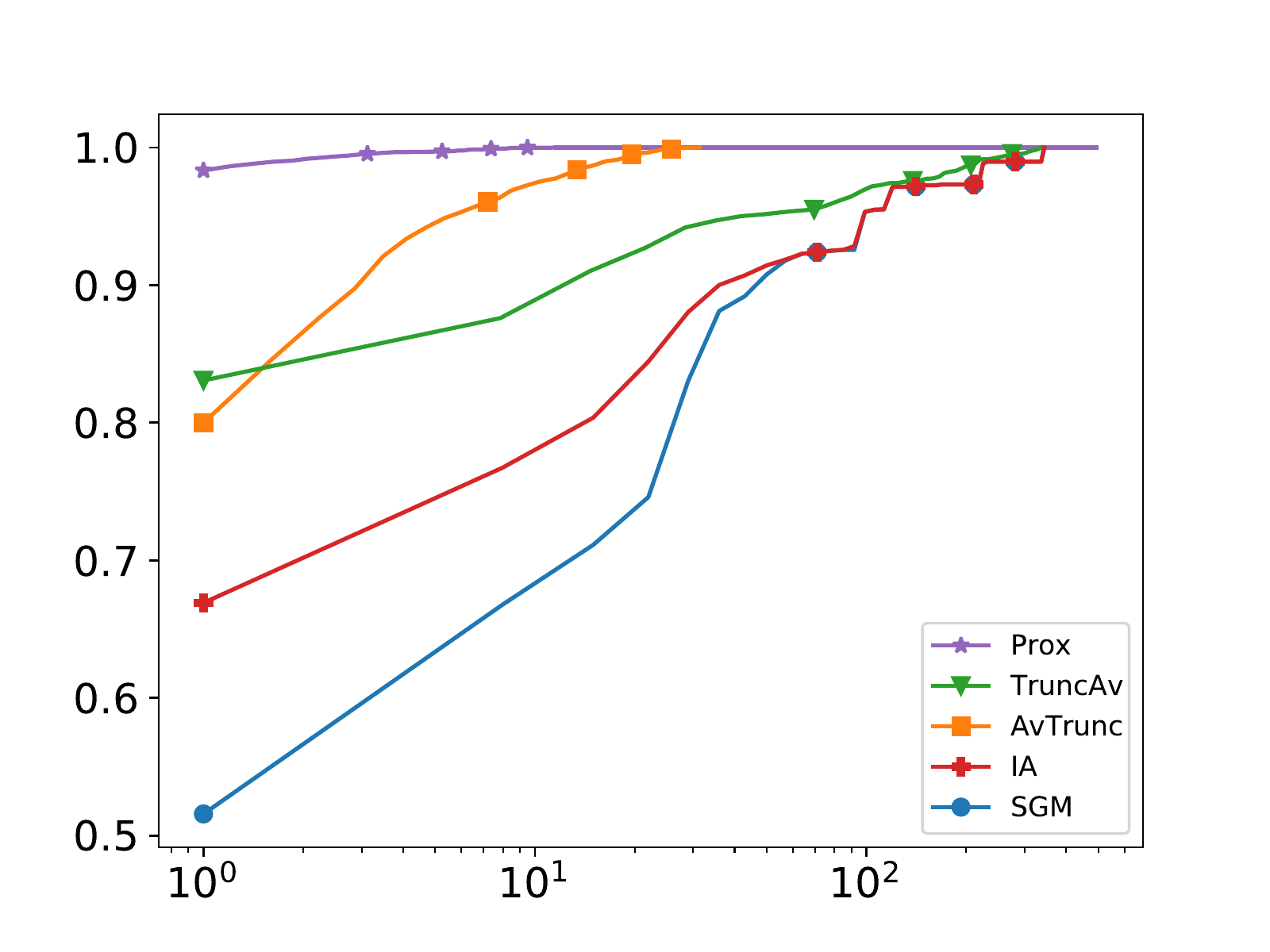}
      		\put(-2,8){
          \rotatebox{90}{{\small Fraction of Experiments}}}
        \put(28,-1){{\small Performance ratio $r$}}
        \put(32,68){{Non-accelerated}}
      \end{overpic} &
      \begin{overpic}[width=.45\columnwidth]{%,grid]{%
      		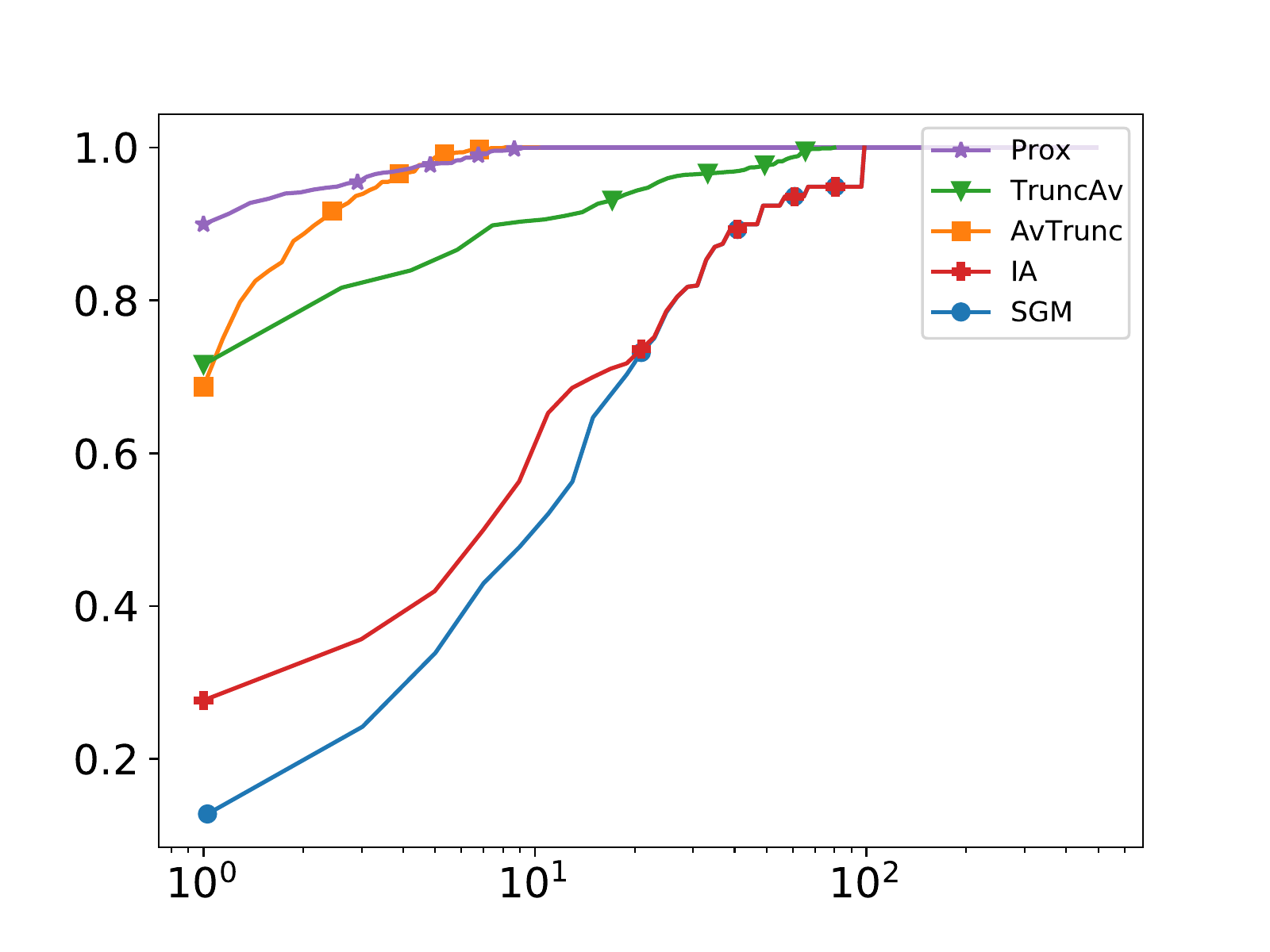}
      		\put(-2,8){
          \rotatebox{90}{{\small Fraction of Experiments}}}
          \put(28,-1){{\small Performance ratio $r$}}
          \put(37,68){{Accelerated}}
      \end{overpic}
    \end{tabular}
  \caption{
    \label{fig:abs-perf-plot}
    Performance profiles for absolute regression. 
    %for (a) noiseless absolute regression, (b) noiseless logistic regression, (c) noisy absolute regression and (d) noisy logistic regression.
  }
  \end{center}
  \vspace{-.5cm}
\end{figure}

\begin{figure}[ht]
  \begin{center}
    \begin{tabular}{cccc}
      \begin{overpic}[width=.45\columnwidth]{%,grid]{%
      		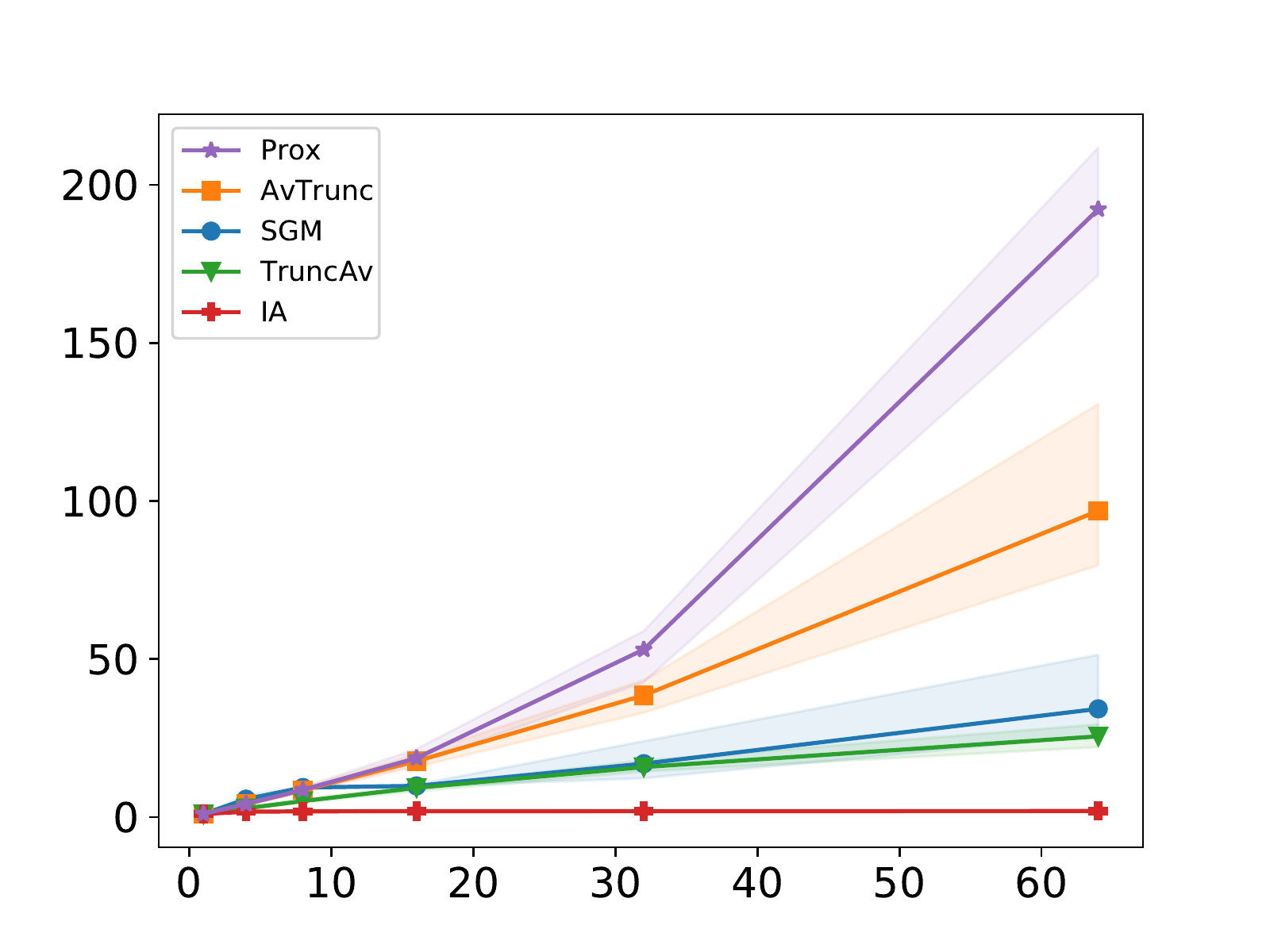}
      		\put(-5,30){
          \rotatebox{90}{{\small Speedup}}}
        \put(28,-1){{\small Minibatch size $\mb$}}
        \put(32,68){{Non-accelerated}}
      \end{overpic} &
      \begin{overpic}[width=.45\columnwidth]{%,grid]{%
      	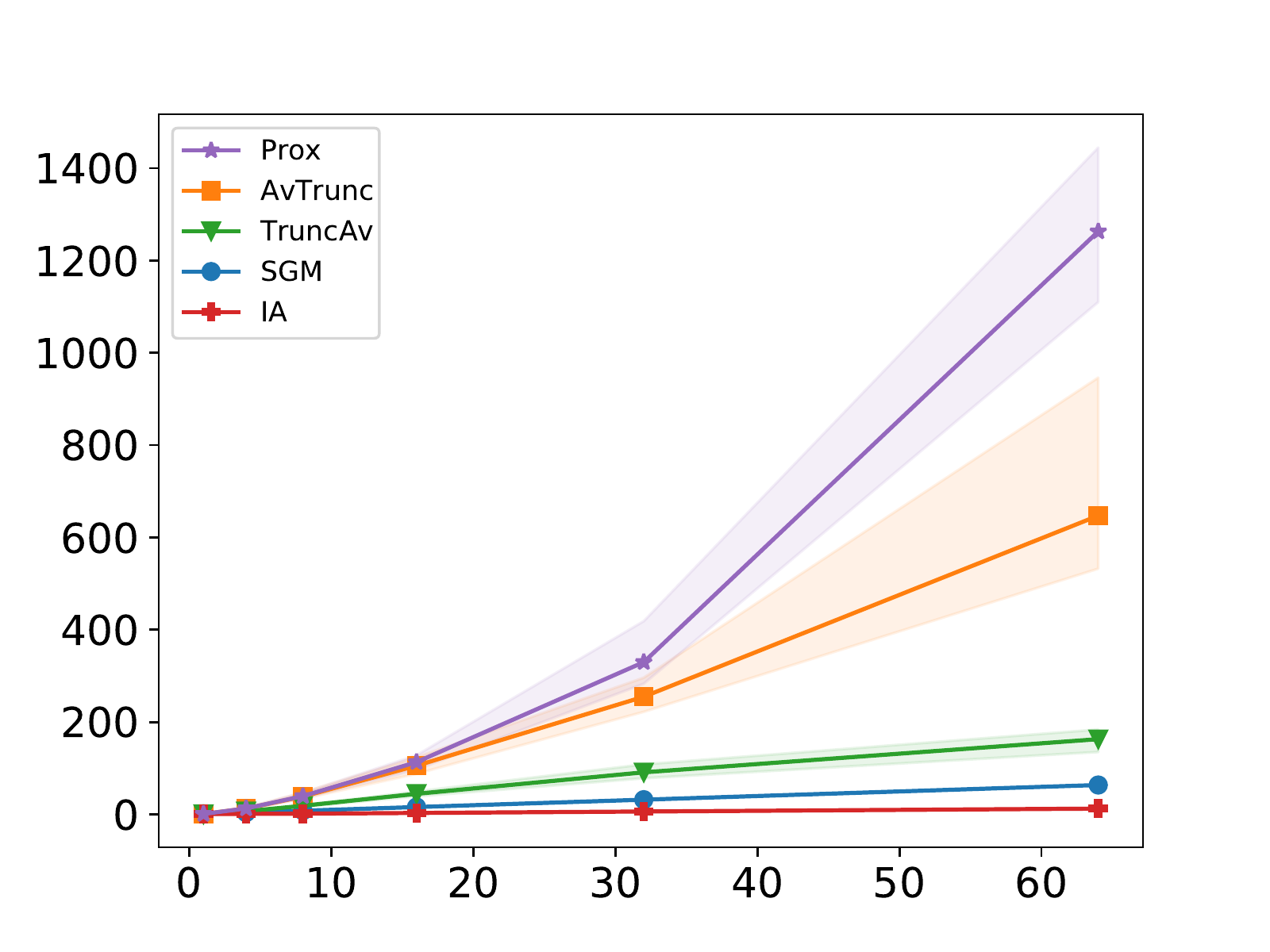}
      		\put(-5,30){
          \rotatebox{90}{{\small Speedup}}}
          \put(28,-1){{\small Minibatch size $\mb$}}
          \put(37,68){{Accelerated}}
      \end{overpic}
    \end{tabular}
  \caption{
  \label{fig:abs-speedup}
   Speed ups with best possible stepsizes  vs.\ batch size for noiseless absolute regression.
    %for (a) noiseless absolute regression, (b) noiseless logistic regression, (c) noisy absolute regression and (d) noisy logistic regression.
  }
  \end{center}
  \vspace{-.5cm}
\end{figure}
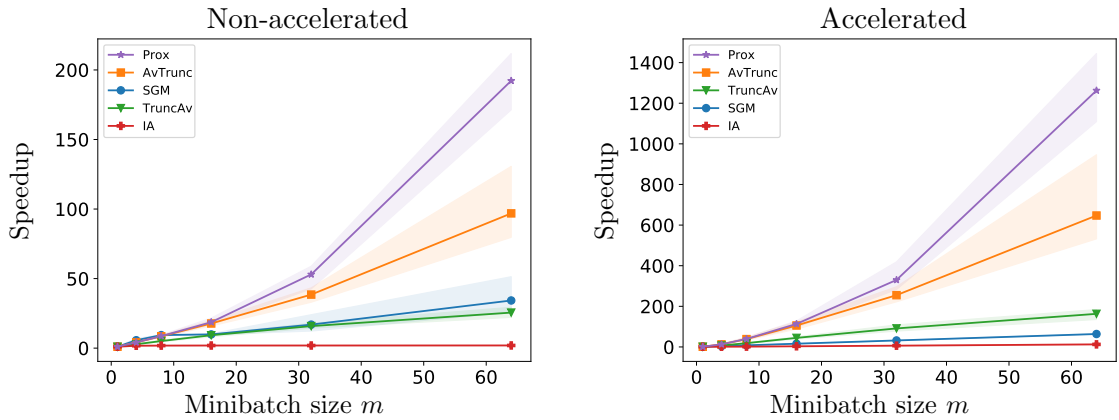

We have $f(x) = \frac{1}{2N}\lone{Ax - b}$. Again we generate rows of $A$
and $x\opt$ i.i.d.\ $\normal(0,I_n)$, setting $b = Ax\opt + \sigma v$ and
drawing $v \sim \laplace(1)^N$. In the noisy setting for our experiments,
we set $\sigma = 0.5$. We provide performance profiles for the
non-accelerated and accelerated algorithms in
\cref{fig:abs-perf-plot}. Similar to the linear regression setting, we see
that \PAM, \PMA, and full-prox, outperform \PIA and SGM. In
\cref{fig:abs-speedup}, we plot the speedup up of each algorithm (relative
to minibatch size $\mb = 1$) against
minibatch size in the noiseless setting. Here, we see the linear improvement
in convergence rate our theoretical results predict,
but there is a superlinear region for large minibatches $\mb > 32$; while
our theory does not predict this, this is because once $\mb \ge n$, a single
step of the stochastic proximal point method can perfectly solve the problem.

\subsection{Logistic Regression} 
\begin{figure}[ht]
  \begin{center}
    \begin{tabular}{cccc}
      \begin{overpic}[width=.45\columnwidth]{%,grid]{%
      		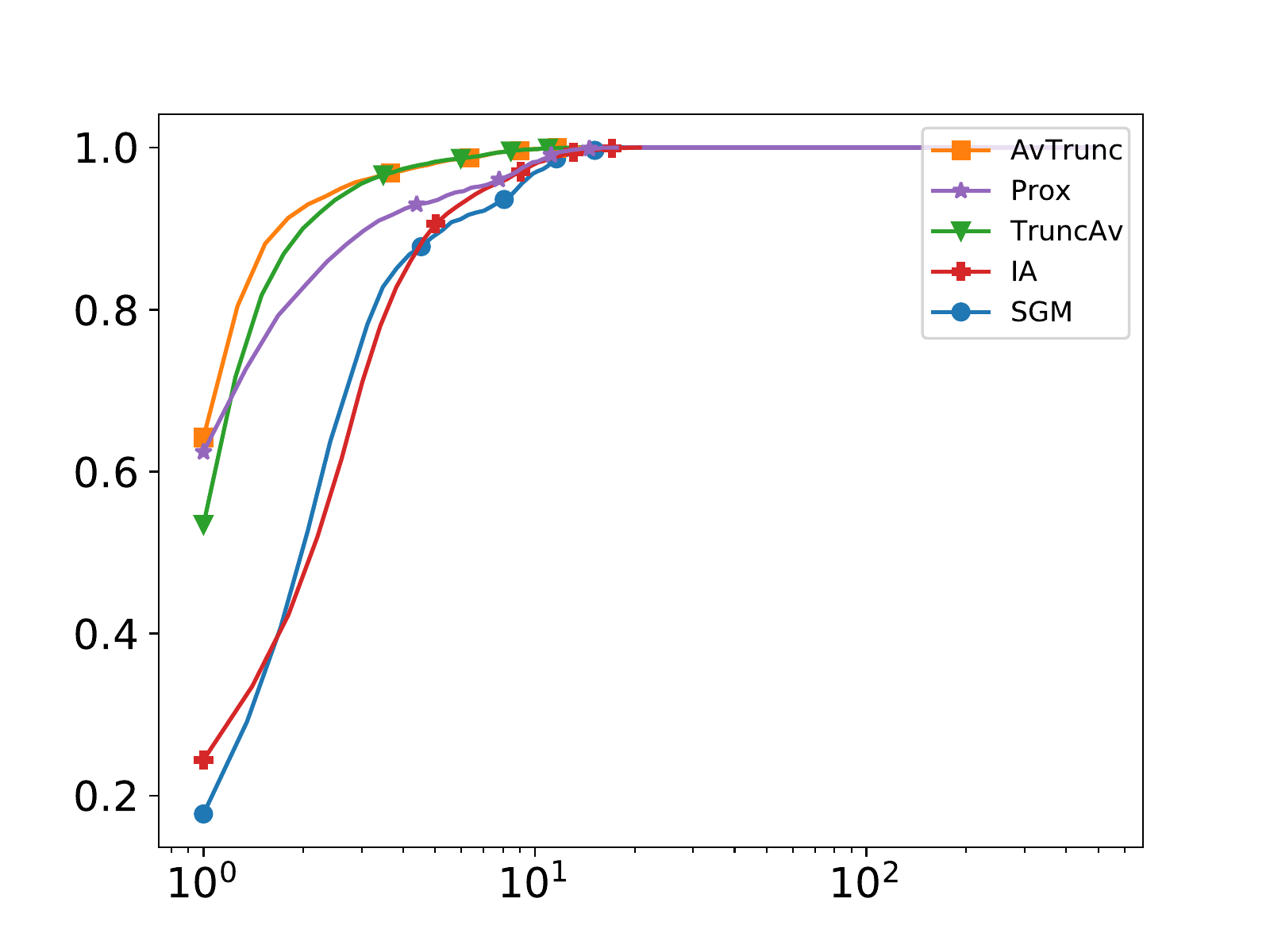}
      		\put(-2,8){
          \rotatebox{90}{{\small Fraction of Experiments}}}
        \put(28,-1){{\small Performance ratio $r$}}
        \put(32,68){{Non-accelerated}}
      \end{overpic} &
      \begin{overpic}[width=.45\columnwidth]{%,grid]{%
      		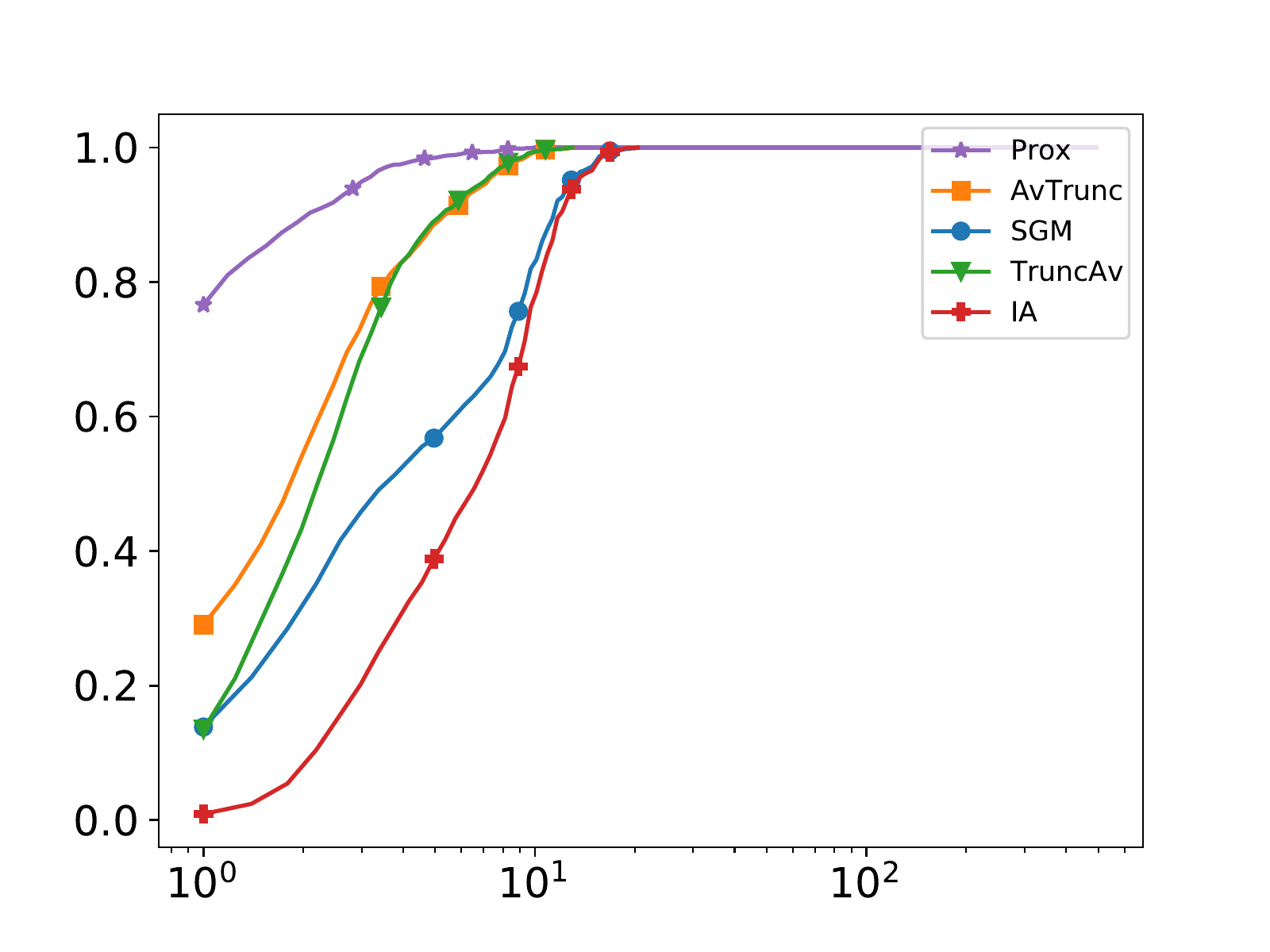}
      		\put(-2,8){
          \rotatebox{90}{{\small Fraction of Experiments}}}
          \put(28,-1){{\small Performance ratio $r$}}
          \put(37,68){{Accelerated}}
      \end{overpic}
    \end{tabular}
  \caption{
    \label{fig:log-perf-plot}
   Performance profiles for logistic regression. 
    %for (a) noiseless absolute regression, (b) noiseless logistic regression, (c) noisy absolute regression and (d) noisy logistic regression.
  }
  \end{center}
  \vspace{-.5cm}
\end{figure}

We have $f(x) = \frac{1}{2N}\sum_{i = 1}^{N}\log(1 + \exp(-b_i \<a_i, x
\>))$. We generate rows of $A$ and $x\opt$ i.i.d.\ $\normal(0,I_n)$,
setting $b_i = \sign(\<a_i, x\opt\>)$. To add noise, we flip each label
$b_i$ independently with probability $p = .01$. We again plot
performance profiles in \cref{fig:log-perf-plot}. The fully proximal,
\PAM, and \PMA methods are noticeably more robust than \PIA and SGM. In the
non-accelerated case, \PAM and \PMA even outperform the fully proximal
method, whereas the fully proximal method outperforms \PAM and \PMA in the
accelerated case. This performance boost from acceleration especially
for the stochastic proximal-point methods may be worthy of further
investigation.

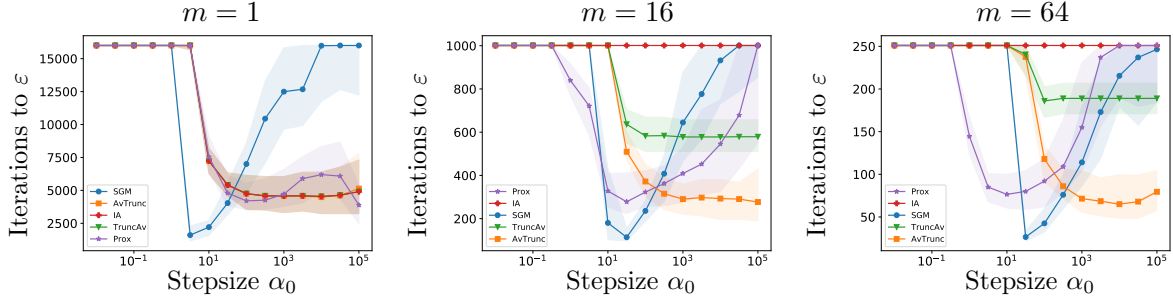
\begin{figure}[ht]
  \begin{center}
    \begin{tabular}{ccc}
      \begin{overpic}[width=.3\columnwidth]{%,grid]{%
      		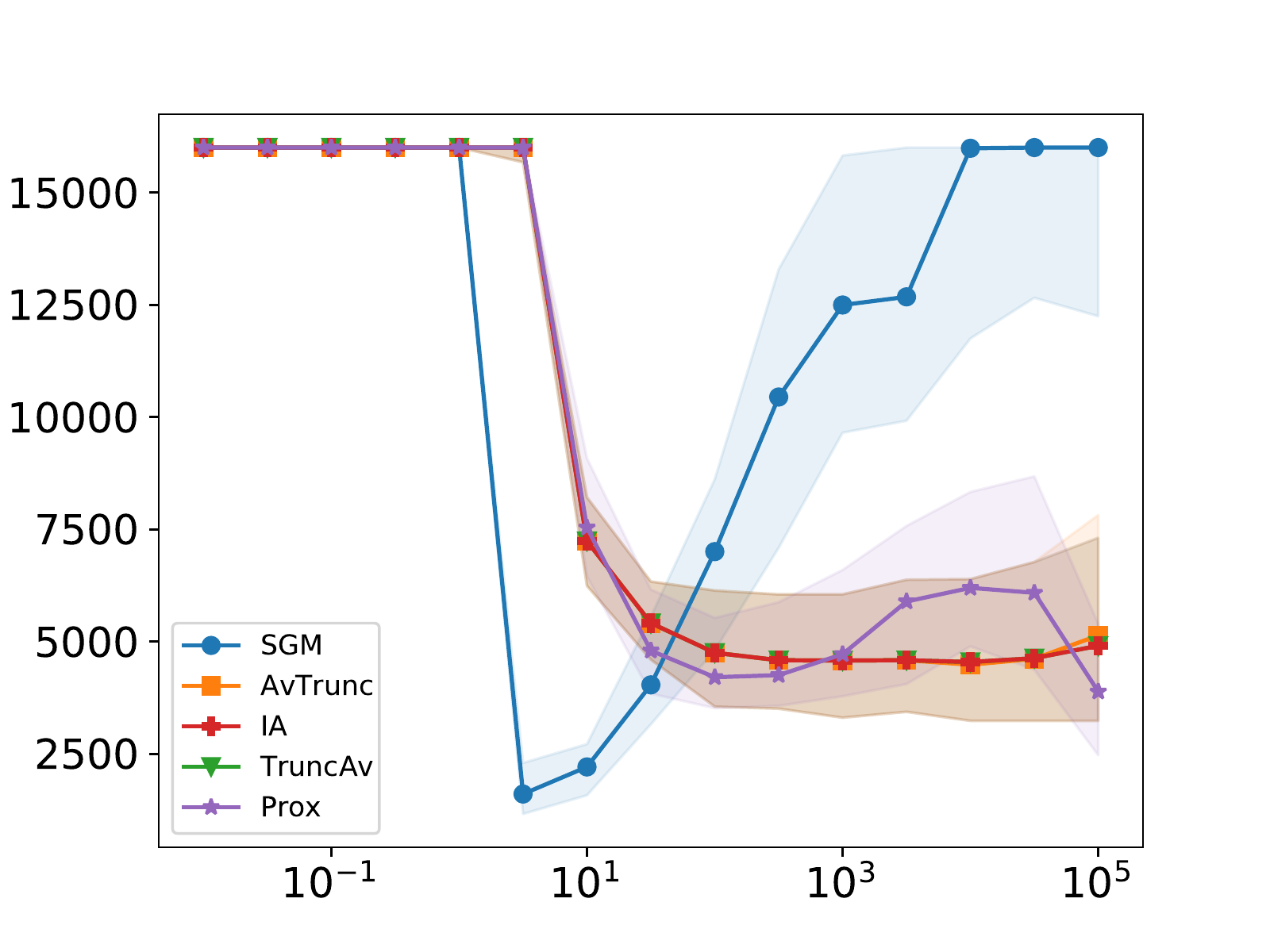}
      		\put(-10,12){
          \rotatebox{90}{{\small Iterations to $\varepsilon$}}}
      		\put(35,-2){{\small Stepsize $\stepsize_0$}}
      		\put(39,70){{$m=1$}}
      \end{overpic} &
      \begin{overpic}[width=.3\columnwidth]{%,grid]{%
      		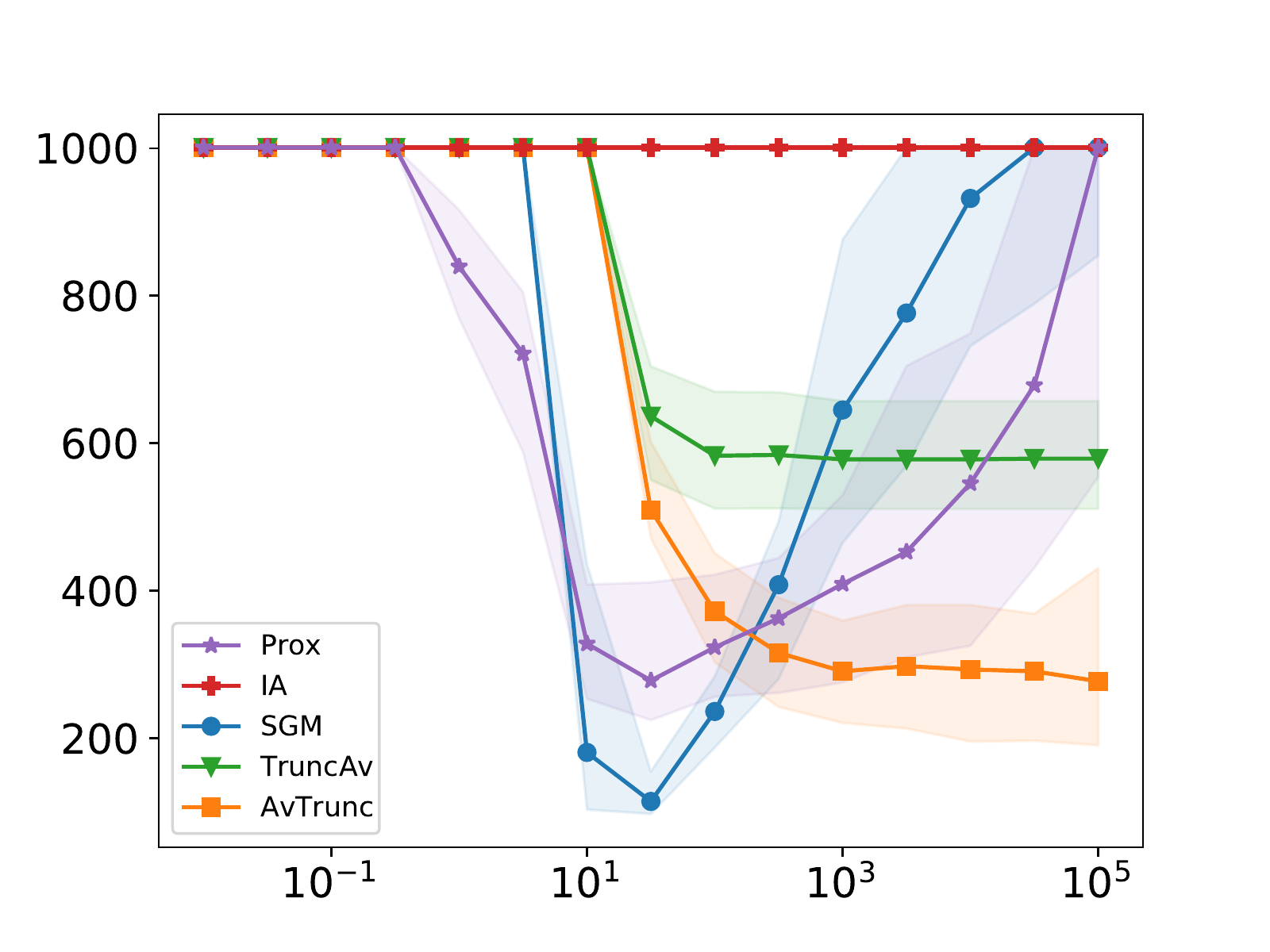}
      		\put(-10,12){
          \rotatebox{90}{{\small Iterations to $\varepsilon$}}}
      		\put(35,-2){{\small Stepsize $\stepsize_0$}}
      		\put(39,70){{$m=16$}}
      \end{overpic} &
      \begin{overpic}[width=.3\columnwidth]{%,grid]{%
      		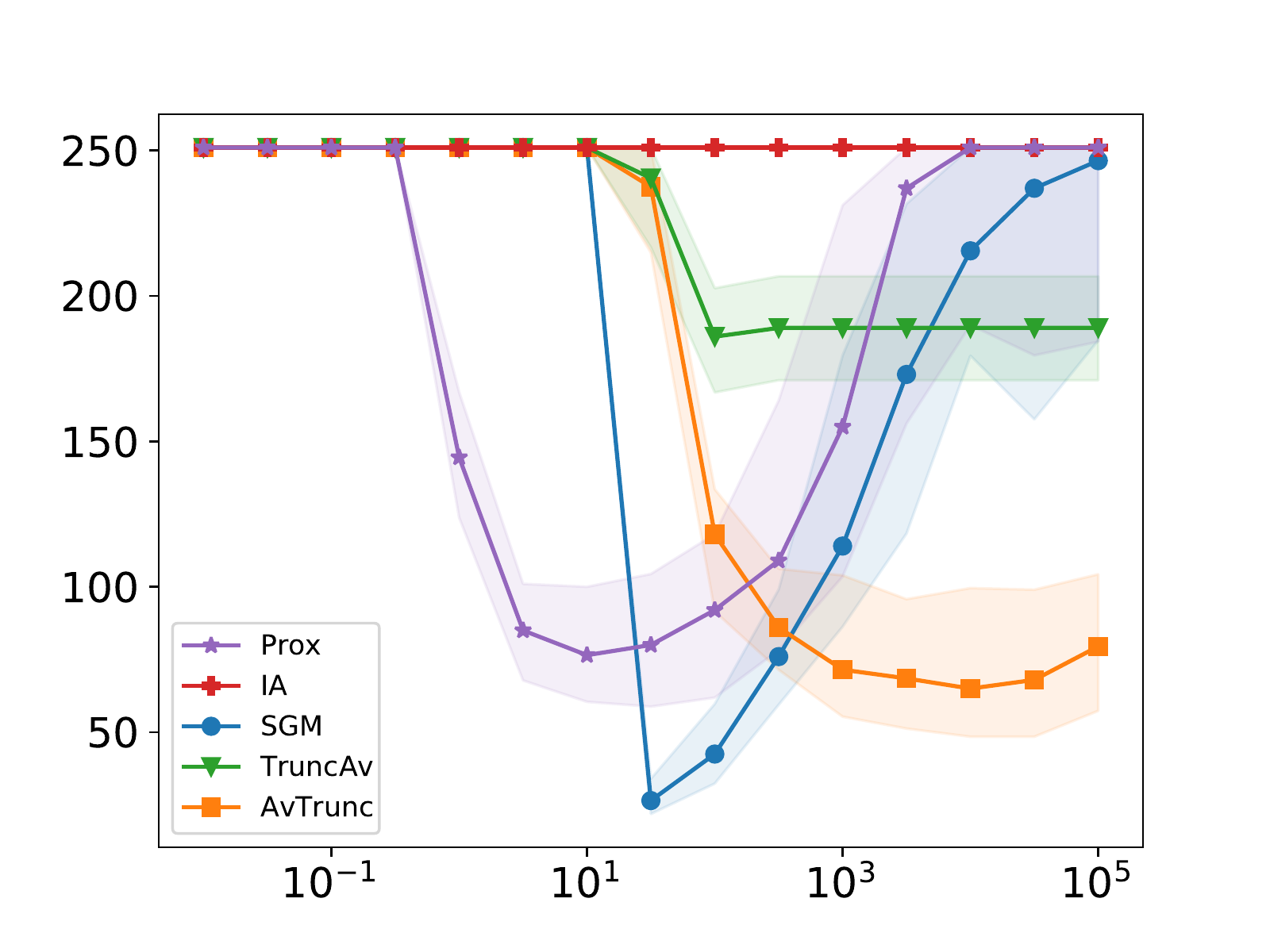}
      		\put(-10,12){
          \rotatebox{90}{{\small Iterations to $\varepsilon$}}}
      		\put(35,-2){{\small Stepsize $\stepsize_0$}}
      		\put(38,70){{$m=64$}}
      \end{overpic}
    \end{tabular}
  \caption{
    \label{fig:noisy-log-iteration}
    Time to convergence of the non-accelerated methods
    vs.\ initial stepsizes 
    for noiseless logistic regression 
    %. (a) $\mb = 1$, (b) $\mb = 8$, (c) $\mb = 32$
    }
  \end{center}
  \vspace{-.3cm}
\end{figure}

\appendix

\section{Proofs of non-asymptotic upper bounds}

We collect our proofs of Theorem~\ref{theorem:non-asymptotic}
and~\ref{theorem:accelerated} in this section.
Both rely on a standard
claim on minimizers of sums of convex functions,
which we state and prove here for convenience.

\begin{claim}
  \label{claim:one-step-minimizer}
  Let $u$ and $\psi$ be convex,  $\psi$ be differentiable
  on $\xdomain$, and
  $D_\psi(x, y) = \psi(x) - \psi(y) - \<\nabla \psi(y), x - y\>$.
  If $x^+$ minimizes $u(x) + \psi(x)$ over $x \in \xdomain$, then
  \begin{equation*}
    u(x^+) + \psi(x^+) \le u(x) + \psi(x) - D_\psi(x, x^+)
    ~~ \mbox{for~all~} x \in \xdomain.
  \end{equation*}
\end{claim}
\begin{proof}
  By convexity and the optimality of $x^+$, there
  exists $u'(x^+) \in \partial u(x^+)$
  such that
  $\<u'(x^+) + \grad \psi(x^+), x - x^+\> \ge 0$ for all $x \in \xdomain$.
  Using the standard first-order convexity inequality, we thus
  obtain
  \begin{align*}
    u(x) & \ge
    u(x^+) + \<u'(x^+), x - x^+\> \\
    & = u(x^+) + \<u'(x^+) + \grad\psi(x^+), x - x^+\>
    - \<\grad\psi(x^+), x - x^+\> \\
    & \ge u(x^+) - \<\grad \psi(x^+), x - x^+\> \\
    & = u(x^+) + \psi(x^+) - \psi(x) + D_\psi(x, x^+),
  \end{align*}
  as desired.
\end{proof}

\subsection{Proof of Theorem~\ref{theorem:non-asymptotic}}
\label{sec:proof-non-asymptotic}

The key to the proof, as is familiar from other analyses of
such methods~\cite{Zinkevich03, NemirovskiJuLaSh09, Lan12, DavisDr19,
  AsiDu19siopt}, is a one-step progress bound.
\begin{lemma}
  \label{lemma:smooth-progress}
  Let the conditions of Theorem~\ref{theorem:non-asymptotic} hold, and
  define
  the function value errors $e_k = [F(x\opt;\statrv_k) - f(x\opt)] -
  [F(x_k; \statrv_k) - f(x_k)]$. Then
  \begin{align*}
    \lefteqn{f(x_{k+1}) - f(x\opt)} \\
    & \le \frac{1}{\stepsize_k} \left[\divergence(x\opt, x_k)
      - \divergence(x\opt, x_{k + 1}) \right]
    + \fnerror_k  
    + \frac{1}{2\eta_k}
    \dnorm{\grad F(x_k;\statrv_k) - \grad f(x_k)}^2. 
  \end{align*}
\end{lemma}
\begin{proof}
  Setting $u(\cdot) = F_{x_k}(\cdot; \statrv_k)$ and $\psi(x) =
  \frac{1}{\stepsize_k} \divergence(x, x_k)$ in
  Claim~\ref{claim:one-step-minimizer}, and taking $x^+ = x_{k + 1}$ and
  $x = x\opt$, we have the progress bound
  \begin{equation}
    \label{eqn:will-need-subs}
    F_{x_k}(x_{k + 1}; \statrv_k)
    +
    \frac{1}{\stepsize_k}
    \divergence(x_{k + 1}, x_k)
    \le F_{x_k}(x\opt; \statrv_k)
    + \frac{1}{\stepsize_k} \left[\divergence(x\opt, x_k)
      - \divergence(x\opt, x_{k + 1})\right].
  \end{equation}
  We turn to bounding the difference $F_{x_k}(x\opt; \statrv_k)
  - F_{x_k}(x_{k + 1}; \statrv_k)$.
  Let $g_k = \grad F(x_k; \statrv_k)$
  and define the gradient error
  $\graderror_k \defeq g_k - \nabla f(x_k)$.
  Using the convexity of $F_{x_k}(\cdot;\statrv_k)$ and
  recalling that $g_k \in \partial F_{x_k}(x_k; \statrv_k)$
  as in our discussion following Condition~\ref{cond:lower-model},
  we have $F_{x_k}(x_{k+1}; \statrv_k) \geq F_{x_k}(x_{k}; \statrv_k)
  + \< g_k, x_{k+1} - x_k\>$. As a consequence, we have
  \begin{align*}
    \lefteqn{F_{x_k}(x\opt; \statrv_k) - F_{x_k}(x_{k+1}; \statrv_k)
      \le F_{x_k}(x\opt;\statrv_k) - F(x_{k};\statrv_k) +
      \<g_k, x_k - x_{k + 1}\>} \\
    & \qquad\qquad\qquad = F_{x_k}(x\opt; \statrv_k) - F(x_k; \statrv_k)
    + \<\nabla f(x_k), x_k - x_{k + 1} \> 
    + \<\graderror_k, x_k - x_{k + 1}\> \\
    & \qquad\qquad\quad~~\!\stackrel{\textup{\ref{cond:lower-model}}}{\le}
    F(x\opt; \statrv_k) - F(x_k; \statrv_k)
    + \<\nabla f(x_k), x_k - x_{k + 1} \> 
    + \<\graderror_k, x_k - x_{k + 1}\> \\
    & \qquad\qquad\qquad = f(x\opt) - f(x_k)
    + \<\nabla f(x_k), x_k - x_{k + 1} \> 
    + \fnerror_k
    + \<\graderror_k, x_k - x_{k + 1}\>,
  \end{align*}
  where we used the error $\fnerror_k = [F(x\opt;\statrv_k) -
    f(x\opt)] - [F(x_k;\statrv_k) - f(x_k)]$.
  Finally, the smoothness of $f$ implies
  $f(x_{k + 1}) \le f(x_k) + \<\nabla f(x_k), x_{k + 1} - x_k\>
  + \frac{\lipgrad}{2} \norm{x_k - x_{k+1}}^2$, so
  \begin{align*}
    \lefteqn{F_{x_k}(x\opt;\statrv_k) - F_{x_k}(x_{k+1};\statrv_k)} \\
    & \qquad \le f(x\opt) - f(x_{k+1}) + \frac{\lipgrad}{2}\norm{x_k - x_{k+1}}^2 +
    \fnerror_k + \langle \graderror_k,x_{k} - x_{k+1}\rangle.
  \end{align*}
  Substituting this into inequality \eqref{eqn:will-need-subs} and
  rearranging, we obtain
  \begin{equation}
    \label{eqn:intermediate-bound}
    \begin{split}
      f(x_{k+1}) - f(x\opt) &\le \frac{1}{\stepsize_k}
      \left[\divergence(x\opt, x_k) - \divergence(x\opt, x_{k + 1})
        - \divergence(x_k, x_{k + 1}) \right] \\
      & \qquad ~ + \fnerror_k + \< \graderror_k,x_{k} - x_{k+1}\>
      + \frac{\lipgrad}{2}\norm{x_k - x_{k + 1}}^2.
    \end{split}
  \end{equation}

  We apply the Fenchel-Young inequality to control the error
  $\<\graderror_k, x_k - x_{k + 1}\>$:
  we have $\<\graderror_k, x_k - x_{k + 1}\>
  \le \frac{1}{2 \eta_k} \dnorm{\graderror_k}^2
  + \frac{\eta_k}{2} \norm{x_k - x_{k + 1}}^2$, so
  \begin{align*}
    f(x_{k+1}) - f(x\opt) &
    \le \frac{1}{\stepsize_k}\left[\divergence(x\opt, x_k)
      - \divergence(x\opt, x_{k+1}) \right] \\
    & \qquad +  \fnerror_k + \frac{1}{2\eta_k} \dnorm{\graderror_k}^2
    + \frac{\lipgrad + \eta_k}{2} \norm{x_k - x_{k+1}}^2
    - \frac{1}{\stepsize_k} \divergence(x_k, x_{k+1}),
  \end{align*}
  which with $\stepsize_k = \frac{1}{\lipgrad + \eta_k}$
  gives the lemma once we apply the strong convexity
  of $\distgen$, that is, that
  $\divergence(x_k, x_{k + 1}) \ge \half \norm{x_k - x_{k + 1}}^2$.
\end{proof}

To complete the proof of the theorem,
we simply sum Lemma~\ref{lemma:smooth-progress}:
\begin{align*}
  \sum_{i = 1}^{k}[f(x_{i+1}) - f(x\opt)]
  &\le \sum_{i = 2}^{k}\left(\frac{1}{\stepsize_i} - \frac{1}{\stepsize_{i-1}}\right) \divergence(x\opt, x_i)
  - \frac{1}{\stepsize_{k + 1}} \divergence(x\opt, x_{k+1}) \\
  &\qquad + \frac{1}{\stepsize_1} \divergence(x\opt, x_1)
  + \sum_{i = 1}^{k} \fnerror_i   + \sum_{i = 1}^{k}\frac{1}{2\eta_i} 
  \dnorm{\grad F(x_i;\statrv_i) - \grad f(x_i)}^2\\ 
  & \le \frac{\radius^2}{\stepsize_k} + \sum_{i = 1}^{k}
  \fnerror_i   + \sum_{i = 1}^{k}\frac{1}{2\eta_i} 
  \dnorm{\grad F(x_i;\statrv_i) - \grad f(x_i)}^2.
\end{align*}
Taking expectations and using that $\E[\fnerror_k] = 0$
and $\stepsize_k = \frac{1}{\lipgrad + \eta_k}$ gives the theorem.

\subsection{Proof of Theorem~\ref{theorem:accelerated}}
\label{sec:proof-accelerated}

The proof is somewhat analogous to that of
Theorem~\ref{theorem:non-asymptotic}, in that we begin with a
deterministic one-step progress bound and then iterate the bound.  In
analogy to Lemma~\ref{lemma:smooth-progress}, we rely on the conditionally
mean-zero function and gradient errors
\begin{equation*}
  \fnerror_k \defeq F(x\opt; \statrv_k) - f(x\opt)
  + f(y_k) - F(y_k; \statrv_k)
  ~~ \mbox{and} ~~
  \graderror_k \defeq \nabla f(y_k) - \nabla F(y_k; \statrv_k).
\end{equation*}
We have the one-step progress bound
\begin{lemma}
  \label{lemma:one-step-accelerated}
  Let $\stepsize_k \le \frac{1}{\lipgrad \theta_k + \eta_k}$
  and $\Delta_k = f(x_k) + r(x_k) - f(x\opt) - r(x\opt)$.
  Then
  \begin{align*}
    \lefteqn{\Delta_{k+1}} \\
    & \le (1 - \theta_k)
    \Delta_k + \theta_k \left[\fnerror_k + \<\graderror_k, z_k - y_k\>
      + \frac{\dnorm{\graderror_k}^2}{2 \eta_k} +
      \frac{1}{\stepsize_k}
      \left(\divergence(x\opt, z_k) - \divergence(x\opt, z_{k + 1}) \right)
      \right].
  \end{align*}
\end{lemma}
\begin{proof}
  We follow the proof of Tseng, Proposition 1~\cite{Tseng08}. For
  shorthand, let
  \begin{equation*}
    \linapprox(x, y) \defeq f(y) + \<\nabla f(y), x - y\> + r(x),
  \end{equation*}
  which linearly approximates $f$ and does not approximate the additive
  component $r$. Then by the $\lipgrad$-smoothness of $\nabla f$, we
  obtain
  \begin{align}
    \nonumber
    f(x_{k + 1}) \, + \,& r(x_{k+1})
    \le \linapprox(x_{k+1}, y_k) + \frac{\lipgrad}{2} \norm{x_{k + 1} - y_k}^2\\
    & = \linapprox((1 - \theta_k) x_k + \theta_k z_{k+1}, y_k)
    + \frac{\lipgrad \theta_k^2}{2}
    \norm{z_{k + 1} - z_k}^2 \nonumber \\
    & \stackrel{(i)}{\le} (1 - \theta_k) \linapprox(x_k, y_k)
    + \theta_k \linapprox(z_{k+1}, y_k)
    + \frac{\lipgrad \theta_k^2}{2}
    \norm{z_{k + 1} - z_k}^2 \nonumber \\
    & \stackrel{(ii)}{\le} \! (1 - \theta_k) (f(x_k) + r(x_k))
    + \theta_k \left[\linapprox(z_{k+1}, y_k)
      + \frac{\lipgrad \theta_k}{2}
      \norm{z_{k + 1} - z_k}^2\right],
    \label{eqn:initial-bound}
  \end{align}
  where the inequality $(i)$ used that $r$ is convex and
  $(ii)$ that $f$ is convex.

  We consider the final two terms in the bound~\eqref{eqn:initial-bound},
  and with function and gradient errors $\fnerrory_k \defeq f(y_k)
  - F(y_k; \statrv_k)$ and $\graderror_k \defeq \nabla f(y_k) - \nabla
  F(y_k; \statrv_k)$, we expand the first in terms of the random samples
  to write
  \begin{align}
    \nonumber
    \linapprox(z_{k+1}, y_k)
    & = F(y_k; \statrv_k) + \<\nabla F(y_k; \statrv_k), z_{k+1}\! - y_k\>
    + r(z_{k + 1})
    + \fnerrory_k \!+\! \<\graderror_k, z_{k + 1} \!- y_k\> \\
    & \le F_{y_k}(z_{k + 1}; \statrv_k) + r(z_{k + 1})
    + \fnerrory_k + \<\graderror_k, z_{k + 1} - y_k\>,
    \label{eqn:linear-bound-with-error}
  \end{align}
  where the inequality uses that the models $F_{y_k}$ necessarily upper
  bound the first-order (linear) approximation to $F$ at $y_k$ (recall the
  discussion following Condition~\ref{cond:lower-model}).  To control
  term~\eqref{eqn:linear-bound-with-error}, we apply
  Claim~\ref{claim:one-step-minimizer} with $u(x) = F_{y_k}(x; \statrv_k)$,
  $\psi(x) = \divergence(x, z_k)$, and $x^+ = z_{k+1}$, so that
  inequality~\eqref{eqn:linear-bound-with-error} implies
  \begin{align*}
    \linapprox(z_{k + 1}, y_k)
    & \le F_{y_k}(x; \statrv_k) + r(x)
    + \frac{1}{\stepsize_k} \left[\divergence(x, z_k) - \divergence(x,
      z_{k + 1}) - \divergence(z_{k + 1}, z_k) \right] \\
    & \qquad ~ + \fnerrory_k + \<\graderror_k, z_{k + 1} - y_k\>
  \end{align*}
  for any $x \in \xdomain$.
  Rearranging terms and using the Fenchel-Young inequality
  to see that
  \begin{equation*}
    \<\graderror_k, z_{k + 1} - y_k\>
    = \<\graderror_k, z_k - y_k\>
    + \<\graderror_k, z_{k + 1} - y_k\>
    \le \<\graderror_k, z_k - y_k\>
    + \frac{\dnorm{\graderror_k}^2}{2 \eta_k}
    + \frac{\eta_k}{2} \norm{z_{k + 1} - z_k}^2
  \end{equation*}
  and using the strong convexity bound
  $\divergence(z_{k + 1}, z_k) \ge \half \norm{z_{k + 1} - z_k}^2$
  then implies
  \begin{align*}
    \linapprox(z_{k + 1}, y_k)
    & \le 
    F_{y_k}(x; \statrv_k) + r(x)
    + \frac{1}{\stepsize_k} \left[\divergence(x, z_k) - \divergence(x,
      z_{k + 1})\right]
    + \fnerrory_k \\
    & \qquad + \<\graderror_k, z_k - y_k\>
    + \frac{\dnorm{\graderror_k}^2}{2 \eta_k}
    + \frac{\eta_k}{2}
    \norm{z_{k + 1} - z_k}^2
    - \frac{1}{2 \stepsize_k}
    \norm{z_{k + 1} - z_k}^2.
  \end{align*}
  Our modeling assumptions guarantee that $F_{y_k}(x; \statrv_k)
  \le F(x; \statrv_k)$, so writing the
  function error $\fnerror_k = F(x; \statrv_k) - f(x)
  + f(y_k) - F(y_k; \statrv_k)$ and
  substituting this upper bound on $\linapprox(z_{k + 1}, y_k)$
  into the bound~\eqref{eqn:initial-bound}
  gives the single-step progress guarantee
  \begin{align*}
    f(x_{k + 1}) + r(x_{k + 1})
    & \le (1 - \theta_k) (f(x_k) + r(x_k))
    + \theta_k (f(x) + r(x)) \\
    & ~~~ + \theta_k \bigg[\fnerror_k
      + \<\graderror_k, z_k - y_k\>
      + \frac{\dnorm{\graderror_k}^2}{2 \eta_k}
      + \frac{1}{\stepsize_k} \left[\divergence(x, z_k) -
        \divergence(x, z_{k + 1}) \right] \\
      & \qquad ~~~~
      + \frac{\lipgrad \theta_k + \eta_k}{2}
      \norm{z_{k + 1} - z_k}^2
      - \frac{1}{\stepsize_k}
      \norm{z_{k + 1} - z_k}^2 \bigg].
  \end{align*}
  Any stepsize
  $\stepsize_k \le \frac{1}{L \theta_k + \eta_k}$ cancels the
  the $\norm{z_{k + 1} - z_k}^2$ terms, and
  setting $x = x\opt$ gives the lemma.
\end{proof}

Iterating Lemma~\ref{lemma:one-step-accelerated} with $\Delta_k = f(x_k) +
r(x_k) - f(x\opt) - r(x\opt) \ge 0$ yields the following deterministic
convergence guarantee.
\begin{lemma}
  \label{lemma:convergence-with-noise}
  Let the conditions of Theorem~\ref{theorem:accelerated} hold.
  Define
  the error terms
  $\lotserror_k \defeq
  \fnerror_k + \<\graderror_k, z_k - y_k\>
  + \frac{\dnorm{\graderror_k}^2 - \gradvariance^2}{2 \stepsize_k}$.
  Then
  \begin{equation*}
    \frac{1}{\theta_k^2}
    \left[f(x_{k + 1}) + r(x_{k + 1})
      - f(x\opt) - r(x\opt)\right]
    \le \sum_{i = 0}^k
    \frac{\gradvariance^2}{2 \theta_i \eta_i}
    + \left(\lipgrad + \frac{\eta_k}{\theta_k}\right) \radius^2
    + \sum_{i = 0}^k \frac{1}{\theta_i} \lotserror_i.
  \end{equation*}
\end{lemma}
\begin{proof}
  Lemma~\ref{lemma:one-step-accelerated} yields
  \begin{align*}
    \frac{1}{\theta_k^2}
    \Delta_{k + 1} & \le
    \frac{1 - \theta_k}{\theta_k^2} \Delta_k
    + \frac{1}{\theta_k \stepsize_k}
    \left[\divergence(x\opt, z_k) - \divergence(x\opt, z_{k+1}) \right]
    + \frac{\gradvariance^2}{2 \eta_k \theta_k} \\
    & \qquad ~ + \frac{1}{\theta_k} \bigg[\underbrace{\fnerror_k
        + \<\graderror_k, z_k - y_k\>
        + \frac{\dnorm{\graderror_k}^2 - \gradvariance^2}{2
          \stepsize_k}}_{\eqdef \lotserror_k} \bigg] \\
    & \le \frac{1}{\theta_{k-1}^2} \Delta_k
    + \frac{1}{\theta_k \stepsize_k}
    \left[\divergence(x\opt, z_k) - \divergence(x\opt, z_{k+1}) \right]
    + \frac{\gradvariance^2}{2 \eta_k \theta_k}
    + \frac{1}{\theta_k} \lotserror_k.
  \end{align*}
  where we recalled that $(1 - \theta_k) / \theta_k^2 \le 1 / \theta_{k-1}^2$.
  Iterating the inequality and using $\frac{1 - \theta_0}{\theta_0^2} = 0$,
  we find that
  \begin{equation*}
    \frac{1}{\theta_k^2} \Delta_{k + 1}
    \le \sum_{i = 0}^k \frac{\gradvariance^2}{2 \eta_i \theta_i}
    + \sum_{i = 0}^k \frac{1}{\theta_i \stepsize_i}
    \left(\divergence(x\opt, z_i) - \divergence(x\opt, z_{i + 1})\right)
    + \sum_{i = 0}^k \frac{1}{\theta_i} \lotserror_i.
  \end{equation*}
  Rearranging the middle summation above as in
  the proof of Theorem~\ref{theorem:non-asymptotic} and
  noting that $\frac{1}{\theta_i \stepsize_i}
  = \lipgrad + \frac{\eta_i}{\theta_i}$ gives
  $\sum_{i = 0}^k \frac{1}{\theta_i \stepsize_i}
  (\divergence(x\opt, z_i) - \divergence(x\opt, z_{i + 1})
  \le \lipgrad \radius^2 + \frac{\eta_k}{\theta_k} \radius^2$,
  as desired.
  %% \begin{align*}
  %%   \lefteqn{\sum_{i = 0}^k \frac{1}{\theta_i \stepsize_i}
  %%     \left(\divergence(x\opt, z_i) - \divergence(x\opt, z_{i + 1})\right)} \\
  %%   & = \sum_{i = 1}^k
  %%   \left[\frac{1}{\theta_i \stepsize_i}
  %%     - \frac{1}{\theta_{i - 1} \stepsize_{i - 1}} \right] \divergence(x\opt,
  %%   z_i)
  %%   + \frac{1}{\theta_0 \stepsize_0} \divergence(x\opt, z_0)
  %%   - \frac{1}{\theta_k \stepsize_k} \divergence(x\opt, z_{k + 1}) \\
  %%   & \le
  %%   \sum_{i = 1}^k \left[\frac{\eta_i}{\theta_i} - \frac{\eta_{i-1}}{\theta_{i-1}}
  %%     \right] \radius^2 + (\lipgrad + \eta_0) \radius^2
  %%   = \lipgrad \radius^2 + \frac{\eta_k}{\theta_k} \radius^2.   
  %% \end{align*}
  %% This gives the lemma.
\end{proof}

Now take
expectations in Lemma~\ref{lemma:convergence-with-noise}. We have
$\E[\lotserror_k] \le 0$, and
\begin{align*}
  \sum_{i = 0}^k \frac{i + 2}{\sqrt{i + 1}}
  & \le \sum_{i = 1}^{k + 1}
  \sqrt{i} + \sum_{i = 1}^{k + 1} \frac{1}{\sqrt{i}} \\
  & \le \int_1^{k + 2} \sqrt{t} dt
  + \int_0^{k + 1} \frac{1}{\sqrt{t}} dt
  =
  \frac{2}{3} ((k + 2)^{3/2} - 1)
  + 2 \sqrt{k + 1}
  \stackrel{(i)}{\le} (k + 2)^{3/2},
\end{align*}
where inequality~$(i)$ holds for $k > 2$.
Multiplying by $\theta_k^2 = 4 / (k + 2)^2$ and using
$\eta_k \theta_k = \eta_0 \frac{2 \sqrt{k}}{k + 2}
\le 2 \eta_0 / \sqrt{k}$
gives the deterministic bound
\begin{equation*}
  \theta_k^2 \sum_{i = 0}^k \frac{\sigma^2}{2 \theta_i \eta_i}
  + \theta_k^2 \left(\lipgrad + \frac{\eta_k}{\theta_k}\right) \radius^2
  \le
  \frac{4 \lipgrad \radius^2}{(k + 2)^2}
  +\frac{2 \radius^2 \eta_0}{\sqrt{k}}
  + \frac{2 \sigma^2}{\eta_0 \sqrt{k + 2}},
\end{equation*}
as desired.

\bibliographystyle{siamplain}
\bibliography{bib}

\begin{thebibliography}{10}

\bibitem{AgarwalBaRaWa12}
{\sc A.~Agarwal, P.~L. Bartlett, P.~Ravikumar, and M.~J. Wainwright}, {\em
  Information-theoretic lower bounds on the oracle complexity of convex
  optimization}, IEEE Transactions on Information Theory, 58 (2012),
  pp.~3235--3249.

\bibitem{AsiDu19}
{\sc H.~Asi and J.~C. Duchi}, {\em The importance of better models in
  stochastic optimization}, Proceedings of the National Academy of Sciences,
  116 (2019), pp.~22924--22930, \url{https://doi.org/10.1073/pnas.1908018116}.

\bibitem{AsiDu19siopt}
{\sc H.~Asi and J.~C. Duchi}, {\em Stochastic (approximate) proximal point
  methods: Convergence, optimality, and adaptivity}, SIAM Journal on
  Optimization, 29 (2019), pp.~2257--2290,
  \url{https://arXiv.org/abs/1810.05633}.

\bibitem{BachMo11}
{\sc F.~Bach and E.~Moulines}, {\em Non-asymptotic analysis of stochastic
  approximation algorithms for machine learning}, in Advances in Neural
  Information Processing Systems 24, 2011, pp.~451--459.

\bibitem{BauschkeBo96}
{\sc H.~Bauschke and J.~Borwein}, {\em On projection algorithms for solving
  convex feasibility problems}, SIAM Review, 38 (1996), pp.~367--426.

\bibitem{BeckTe03}
{\sc A.~Beck and M.~Teboulle}, {\em Mirror descent and nonlinear projected
  subgradient methods for convex optimization}, Operations Research Letters, 31
  (2003), pp.~167--175.

\bibitem{BelkinHsMi18}
{\sc M.~Belkin, D.~Hsu, and P.~Mitra}, {\em Overfitting or perfect fitting?
  {R}isk bounds for classification and regression rules that interpolate}, in
  Advances in Neural Information Processing Systems 31, Curran Associates,
  Inc., 2018, pp.~2300--2311.

\bibitem{BelkinRaTs19}
{\sc M.~Belkin, A.~Rakhlin, and A.~B. Tsybakov}, {\em Does data interpolation
  contradict statistical optimality?}, in Proceedings of the 22nd International
  Conference on Artificial Intelligence and Statistics, 2019, pp.~1611--1619.

\bibitem{Bertsekas11}
{\sc D.~P. Bertsekas}, {\em Incremental proximal methods for large scale convex
  optimization}, Mathematical Programming, Series B, 129 (2011), pp.~163--195.

\bibitem{BottouBo07}
{\sc L.~Bottou and O.~Bousquet}, {\em The tradeoffs of large scale learning},
  in Advances in Neural Information Processing Systems 20, 2007.

\bibitem{CarmonDuHiSi19}
{\sc Y.~Carmon, J.~C. Duchi, O.~Hinder, and A.~Sidford}, {\em Lower bounds for
  finding stationary points {I}}, Mathematical Programming, Series A, to appear
  (2019).

\bibitem{ChaturapruekDuRe15}
{\sc S.~Chaturapruek, J.~C. Duchi, and C.~R\'e}, {\em Asynchronous stochastic
  convex optimization: the noise is in the noise and {SGD} don't care}, in
  Advances in Neural Information Processing Systems 28, 2015.

\bibitem{DavisDr19}
{\sc D.~Davis and D.~Drusvyatskiy}, {\em Stochastic model-based minimization of
  weakly convex functions}, SIAM Journal on Optimization, 29 (2019),
  pp.~207--239.

\bibitem{DefazioBaLa14}
{\sc A.~Defazio, F.~Bach, and S.~Lacoste-Julien}, {\em {SAGA}: A fast
  incremental gradient method with support for non-strongly convex composite
  objectives}, in Advances in Neural Information Processing Systems 27, 2014.

\bibitem{DekelGiShXi12}
{\sc O.~Dekel, R.~Gilad-Bachrach, O.~Shamir, and L.~Xiao}, {\em Optimal
  distributed online prediction using mini-batches}, Journal of Machine
  Learning Research, 13 (2012), pp.~165--202.

\bibitem{DolanMo02}
{\sc E.~D. Dolan and J.~J. Mor{\'{e}}}, {\em Benchmarking optimization software
  with performance profiles}, Mathematical Programming, 91 (2002),
  pp.~201--213.

\bibitem{DuchiBaWa12}
{\sc J.~C. Duchi, P.~L. Bartlett, and M.~J. Wainwright}, {\em Randomized
  smoothing for stochastic optimization}, SIAM Journal on Optimization, 22
  (2012), pp.~674--701.

\bibitem{DuchiHaSi11}
{\sc J.~C. Duchi, E.~Hazan, and Y.~Singer}, {\em Adaptive subgradient methods
  for online learning and stochastic optimization}, Journal of Machine Learning
  Research, 12 (2011), pp.~2121--2159.

\bibitem{DuchiRu18c}
{\sc J.~C. Duchi and F.~Ruan}, {\em Stochastic methods for composite and weakly
  convex optimization problems}, SIAM Journal on Optimization, 28 (2018),
  pp.~3229--3259.

\bibitem{FullerMi11}
{\sc S.~Fuller and L.~Millett}, {\em The Future of Computing Performance: Game
  Over or Next Level?}, National Academies Press, 2011.

\bibitem{KarampatziakisLa11}
{\sc N.~Karampatziakis and J.~Langford}, {\em Online importance weight aware
  updates}, in Proceedings of the 27th Conference on Uncertainty in Artificial
  Intelligence, 2011.

\bibitem{KulisBa10}
{\sc B.~Kulis and P.~Bartlett}, {\em Implicit online learning}, in Proceedings
  of the 27th International Conference on Machine Learning, 2010.

\bibitem{KushnerYi03}
{\sc H.~J. Kushner and G.~Yin}, {\em Stochastic Approximation and Recursive
  Algorithms and Applications}, Springer, second~ed., 2003.

\bibitem{Lan12}
{\sc G.~Lan}, {\em An optimal method for stochastic composite optimization},
  Mathematical Programming, Series A, 133 (2012), pp.~365--397.

\bibitem{LevyDu19}
{\sc D.~Levy and J.~C. Duchi}, {\em Necessary and sufficient geometries for
  gradient methods}, in Advances in Neural Information Processing Systems 32,
  2019, \url{https://arxiv.org/abs/1909.10455}.

\bibitem{LiJaDeRoTa17}
{\sc L.~Li, K.~Jamieson, G.~DeSalvo, A.~Rostamizadeh, and A.~Talwalkar}, {\em
  Hyperband: A novel bandit-based approach to hyperparameter optimization},
  Journal of Machine Learning Research, 18 (2017), pp.~1--52.

\bibitem{LinMaHa18}
{\sc H.~Lin, J.~Mairal, and Z.~Harchaoui}, {\em Catalyst acceleration for
  first-order convex optimization: from theory to practice}, Journal of Machine
  Learning Research, 18 (2018).

\bibitem{MaBaBe18}
{\sc S.~Ma, R.~Bassily, and M.~Belkin}, {\em The power of interpolation:
  Understanding the effectiveness of {SGD} in modern over-parametrized
  learning}, in Proceedings of the 35th International Conference on Machine
  Learning, 2018.

\bibitem{ManiaPaPaReRaJo17}
{\sc H.~Mania, X.~Pan, D.~Papailiopoulos, B.~Recht, K.~Ramchandran, and M.~I.
  Jordan}, {\em Perturbed iterate analysis for asynchronous stochastic
  optimization}, SIAM Journal on Optimization, 27 (2017), pp.~2202–--2229.

\bibitem{NeedellWaSr14}
{\sc D.~Needell, R.~Ward, and N.~Srebro}, {\em Stochastic gradient descent,
  weighted sampling, and the randomized {K}aczmarz algorithm}, in Advances in
  Neural Information Processing Systems 27, 2014, pp.~1017--1025.

\bibitem{NemirovskiJuLaSh09}
{\sc A.~Nemirovski, A.~Juditsky, G.~Lan, and A.~Shapiro}, {\em Robust
  stochastic approximation approach to stochastic programming}, SIAM Journal on
  Optimization, 19 (2009), pp.~1574--1609.

\bibitem{NemirovskiYu83}
{\sc A.~Nemirovski and D.~Yudin}, {\em Problem Complexity and Method Efficiency
  in Optimization}, Wiley, 1983.

\bibitem{Nesterov04}
{\sc Y.~Nesterov}, {\em Introductory Lectures on Convex Optimization}, Kluwer
  Academic Publishers, 2004.

\bibitem{RechtReWrNi11}
{\sc F.~Niu, B.~Recht, C.~R\'{e}, and S.~Wright}, {\em Hogwild: a lock-free
  approach to parallelizing stochastic gradient descent}, in Advances in Neural
  Information Processing Systems 24, 2011.

\bibitem{NiuReReWr11}
{\sc F.~Niu, B.~Recht, C.~Re, and S.~Wright}, {\em Hogwild: a lock-free
  approach to parallelizing stochastic gradient descent}, in Advances in Neural
  Information Processing Systems 24, 2011.

\bibitem{OrabonaPa18}
{\sc F.~Orabona and D.~P{\'a}l}, {\em Scale-free online learning}, Theoretical
  Computer Science, 716 (2018), pp.~50--69.

\bibitem{Polyak87}
{\sc B.~T. Polyak}, {\em Introduction to Optimization}, Optimization Software,
  Inc., 1987.

\bibitem{PolyakJu92}
{\sc B.~T. Polyak and A.~B. Juditsky}, {\em Acceleration of stochastic
  approximation by averaging}, SIAM Journal on Control and Optimization, 30
  (1992), pp.~838--855.

\bibitem{RobbinsMo51}
{\sc H.~Robbins and S.~Monro}, {\em A stochastic approximation method}, Annals
  of Mathematical Statistics, 22 (1951), pp.~400--407.

\bibitem{ScamanBaBuLeMa17}
{\sc K.~Scaman, F.~Bach, S.~Bubeck, Y.~T. Lee, and L.~Massouli\'{e}}, {\em
  Optimal algorithms for smooth and strongly convex distributed optimization in
  networks}, in Proceedings of the 34th International Conference on Machine
  Learning, 2017.

\bibitem{ShalevSiSrCo11}
{\sc S.~Shalev-Shwartz, Y.~Singer, N.~Srebro, and A.~Cotter}, {\em Pegasos:
  primal estimated sub-gradient solver for {SVM}}, Mathematical Programming,
  Series B, 127 (2011), pp.~3--30.

\bibitem{StrohmerVe09}
{\sc T.~Strohmer and R.~Vershynin}, {\em A randomized {K}aczmarz algorithm with
  exponential convergence}, Journal of Fourier Analysis and Applications, 15
  (2009), pp.~262--278.

\bibitem{Tseng08}
{\sc P.~Tseng}, {\em On accelerated proximal gradient methods for
  convex-concave optimization}.
\newblock 2008, \url{http://www.mit.edu/~dimitrib/PTseng/papers/apgm.pdf}.

\bibitem{Wainwright19}
{\sc M.~J. Wainwright}, {\em High-Dimensional Statistics: A Non-Asymptotic
  Viewpoint}, Cambridge University Press, 2019.

\bibitem{Zhang04}
{\sc T.~Zhang}, {\em Solving large scale linear prediction problems using
  stochastic gradient descent algorithms}, in Proceedings of the Twenty-First
  International Conference on Machine Learning, 2004.

\bibitem{Zinkevich03}
{\sc M.~Zinkevich}, {\em Online convex programming and generalized
  infinitesimal gradient ascent}, in Proceedings of the Twentieth International
  Conference on Machine Learning, 2003.

\end{thebibliography}

\end{document}